
\magnification=1200

\def\R{{\bf R}}					
\def\upchi{\raise1.2pt\hbox{$\chi$}}		

\def\norm#1{\left\Vert{#1}\right\Vert}                 	
\def\abs#1{\left\vert{#1}\right\vert}                   	
\def\frac#1#2{{{#1}\over{#2}}}			

\def\today{\ifcase\month\or January\or February\or March\or April\or
May\or June\or July\or August\or September\or October\or November\or
December\fi\space\number\day, \number\year}  

\def\intave#1{\int_{#1}\hbox{\llap{$\raise2.3pt\hbox{\vrule
height.9pt width7pt}\phantom{\scriptstyle{#1}}\mkern-2mu$}}}	

\def\div{{\rm div}}			     	
\def\sqr#1#2{{\vcenter{\hrule height.#2pt\hbox{\vrule
     width.#2pt height#1pt\kern#1pt\vrule width.#2pt}\hrule height.#2pt}}}


\def\bull{\vrule height 1ex width .8ex depth -.0ex}

\def\qed{\enspace\hfill $\bull$}    

 5	
 4		
 3	
 2		
\font\eightpt=cmr8				
 5	
 4	
 3	
 2	

\def\flushitem{\par\hang\ztextindent}
\def\ztextindent#1{\indent\llap{\hbox to\parindent{#1\hfill}}\ignorespaces}

\def\nullset{\hbox{\eightpt \O}}
\newif\ifproofmode                      
\proofmodefalse				

\newif\ifforwardreference		
\forwardreferencefalse			

\newif\ifchapternumbers			
\chapternumbersfalse			

\newif\ifcontinuousnumbering		
\continuousnumberingfalse		

\newif\iffigurechapternumbers		
\figurechapternumbersfalse		

\newif\ifcontinuousfigurenumbering	
\continuousfigurenumberingfalse		

\font\eqsixrm=cmr6			
\def\marginstyle{\eqsixrm}		

\newtoks\chapletter			
\newcount\chapno			
\newcount\eqlabelno			
\newcount\figureno			

\chapno=0
\eqlabelno=0
\figureno=0


\def\chapfolio{\ifnum \chapno>0 \the\chapno \else \the\chapletter \fi}


\def\bumpchapno{\ifnum \chapno>-1 \global \advance \chapno by 1
	\else \global \advance \chapno by -1 \setletter\chapno \fi
	\ifcontinuousnumbering \else \global\eqlabelno=0 \fi
	\ifcontinuousfigurenumbering \else \global\figureno=0 \fi}

%


%

\def\tempsetletter#1{\ifcase-#1 {}\or{} \or\chapletter={A}\or\chapletter={B}
  \or\chapletter={C} \or\chapletter={D} \or\chapletter={E}
  \or\chapletter={F} \or\chapletter={G} \or\chapletter={H}
  \or\chapletter={I} \or\chapletter={J} \or\chapletter={K}
  \or\chapletter={L} \or\chapletter={M} \or\chapletter={N}
  \or\chapletter={O} \or\chapletter={P} \or\chapletter={Q}
  \or\chapletter={R} \or\chapletter={S} \or\chapletter={T}
  \or\chapletter={U} \or\chapletter={V} \or\chapletter={W}
  \or\chapletter={X} \or\chapletter={Y} \or\chapletter={Z}\fi}

%

\def\chapshow#1{\ifnum #1>0 \relax #1%
   \else {\tempsetletter{\number#1}\chapno=#1 \chapfolio} \fi}

%

\def\today{\number\day\space \ifcase\month\or Jan\or Feb\or
   Mar\or Apr\or May\or Jun\or Jul\or Aug\or Sep\or
   Oct\or Nov\or Dec\fi, \number\year}

%

\def\initialeqmacro{\ifproofmode
 \headline{\tenrm \today\hfill \jobname\ --- draft\hfill\folio}
     \hoffset=0truein \immediate\openout2=allcrossreferfile \fi
 \ifforwardreference \input labelfile
     \ifproofmode \immediate\openout1=labelfile \fi \fi}


%

\def\chaplabel#1{\bumpchapno \ifproofmode \ifforwardreference
   \immediate\write1{\noexpand\expandafter\noexpand\def
   \noexpand\csname CHAPLABEL#1\endcsname{\the\chapno}}\fi\fi
   \global\expandafter\edef\csname CHAPLABEL#1\endcsname
   {\the\chapno}\ifproofmode\llap{\hbox{\marginstyle #1\ }}\fi\chapfolio}

%
\def\eqnum{\global\advance\eqlabelno by 1
   \eqno(\ifchapternumbers\chapfolio.\fi\the\eqlabelno)}

\def\eqlabel#1{\global\advance\eqlabelno by 1 \ifproofmode\ifforwardreference
 \immediate\write1{\noexpand\expandafter\noexpand\def
 \noexpand\csname EQLABEL#1\endcsname{\the\chapno.\the\eqlabelno?}}\fi\fi
 \global\expandafter\edef\csname EQLABEL#1\endcsname
 {\the\chapno.\the\eqlabelno?} \eqno(\ifchapternumbers\chapfolio.\fi
 \the\eqlabelno)\ifproofmode\rlap{\hbox{\marginstyle #1}}\fi}

\def\leqlabel#1{\global\advance\eqlabelno by 1 \ifproofmode\ifforwardreference
 \immediate\write1{\noexpand\expandafter\noexpand\def
 \noexpand\csname EQLABEL#1\endcsname{\the\chapno.\the\eqlabelno?}}\fi\fi
 \global\expandafter\edef\csname EQLABEL#1\endcsname
 {\the\chapno.\the\eqlabelno?} \leqno(\ifchapternumbers\chapfolio.\fi
 \the\eqlabelno)\ifproofmode\rlap{\hbox{\marginstyle #1}}\fi}

\def\eqalignnum{\global\advance\eqlabelno by 1
   &(\ifchapternumbers\chapfolio.\fi\the\eqlabelno)}

\def\eqalignlabel#1{\global\advance\eqlabelno by1 \ifproofmode
 \ifforwardreference\immediate\write1{\noexpand\expandafter\noexpand\def
 \noexpand\csname EQLABEL#1\endcsname
     {\the\chapno.\the\eqlabelno?}}\fi\fi
 \global\expandafter\edef\csname EQLABEL#1\endcsname
 {\the\chapno.\the\eqlabelno?}&(\ifchapternumbers\chapfolio.\fi
 \the\eqlabelno)\ifproofmode\rlap{\hbox{\marginstyle #1}}\fi}

\def\eqref#1{(\ifundefined{EQLABEL#1}***\ifproofmode\ifforwardreference)%
   \else \write16{ ***Undefined Equation Reference #1*** }\fi
   \else \write16{ ***Undefined Equation Reference #1*** }\fi
   \else \edef\LABxx{\getlabel{EQLABEL#1}}%
   \def\LAByy{\expandafter\stripchap\LABxx}\ifchapternumbers
   \chapshow{\LAByy}.\expandafter\stripeq\LABxx
   \else\ifnum \number\LAByy=\chapno \relax\expandafter\stripeq\LABxx
   \else\chapshow{\LAByy}.\expandafter\stripeq\LABxx\fi\fi)\fi
   \ifproofmode\write2{Equation #1}\fi}

%

\def\fignum{\global\advance\figureno by 1 \relax
   \iffigurechapternumbers\chapfolio.\fi\the\figureno}\

\def\figlabel#1{\global\advance\figureno by 1\relax
 \ifproofmode\ifforwardreference
 \immediate\write1{\noexpand\expandafter\noexpand\def
 \noexpand\csname FIGLABEL#1\endcsname{\the\chapno.\the\figureno?}}\fi\fi
 \global\expandafter\edef\csname FIGLABEL#1\endcsname
 {\the\chapno.\the\figureno?}\iffigurechapternumbers\chapfolio.\fi
 \ifproofmode$^{\hbox{\marginstyle #1}}$\relax\fi\the\figureno}

\def\figref#1{\ifundefined{FIGLABEL#1}!!!!\ifproofmode\ifforwardreference)%
   \else \write16{ ***Undefined Equation Reference #1*** }\fi
   \else \write16{ ***Undefined Equation Reference #1*** }\fi
   \else \edef\LABxx{\getlabel{FIGLABEL#1}}%
   \def\LAByy{\expandafter\stripchap\LABxx}%
   \iffigurechapternumbers\chapshow{\LAByy}.\expandafter\stripeq\LABxx
   \else\ifnum\number\LAByy=\chapno \relax\expandafter\stripeq\LABxx
   \else\chapshow{\LAByy}.\expandafter\stripeq\LABxx\fi\fi
   \ifproofmode\write2{Figure #1}\fi\fi}

%

%

\def\getlabel#1{\csname#1\endcsname}
\def\ifundefined#1{\expandafter\ifx\csname#1\endcsname\relax}
\def\stripchap#1.#2?{#1}
\def\stripeq#1.#2?{#2}

\figurechapternumberstrue  

\chapternumberstrue        

\def\thmlbl#1{\figlabel{#1}}
\def\thmref{\figref}
\def\eqnlbl#1{\leqlabel{#1}}

\def\eqnref#1{\eqref{#1}}
\def\sectionnumber{\chapno}
\def\theoremnumber{\figureno}
\def\equationnumber{\eqlabelno}

\expandafter \def \csname EQLABELassumption\endcsname {1.1?}
\expandafter \def \csname EQLABELnp.1\endcsname {2.2?}
\expandafter \def \csname EQLABELnp1\endcsname {2.3?}
\expandafter \def \csname EQLABELrep\endcsname {2.4?}
\expandafter \def \csname EQLABELnp2\endcsname {2.5?}
\expandafter \def \csname EQLABELnp3\endcsname {2.6?}
\expandafter \def \csname EQLABELareamin\endcsname {2.7?}
\expandafter \def \csname EQLABELnp10\endcsname {2.8?}
\expandafter \def \csname EQLABELnp11\endcsname {2.9?}
\expandafter \def \csname FIGLABELlipdiff\endcsname {2.1?}
\expandafter \def \csname FIGLABELdiffae\endcsname {2.2?}
\expandafter \def \csname EQLABELdiv\endcsname {2.10?}
\expandafter \def \csname FIGLABELmeancur\endcsname {2.3?}
\expandafter \def \csname EQLABELdivmc\endcsname {2.11?}
\expandafter \def \csname EQLABELproblem\endcsname {3.1?}
\expandafter \def \csname EQLABELexcessrn\endcsname {3.2?}
\expandafter \def \csname FIGLABELbk\endcsname {3.1?}
\expandafter \def \csname EQLABELsamevol\endcsname {3.3?}
\expandafter \def \csname EQLABELmceleqk\endcsname {3.4?}
\expandafter \def \csname FIGLABELberegular\endcsname {3.2?}
\expandafter \def \csname FIGLABELbhregular\endcsname {3.3?}
\expandafter \def \csname EQLABELlip\endcsname {3.5?}
\expandafter \def \csname EQLABELc1\endcsname {3.6?}
\expandafter \def \csname EQLABELsuminh\endcsname {3.7?}
\expandafter \def \csname FIGLABELballinmin\endcsname {3.4?}
\expandafter \def \csname EQLABELintersectionisball\endcsname {3.8?}
\expandafter \def \csname EQLABELesamemeasure\endcsname {3.9?}
\expandafter \def \csname EQLABELintmul\endcsname {3.10?}
\expandafter \def \csname FIGLABELdensity\endcsname {3.5?}
\expandafter \def \csname FIGLABELforms\endcsname {3.6?}
\expandafter \def \csname FIGLABELapproxid\endcsname {3.7?}
\expandafter \def \csname EQLABELsumdeg\endcsname {3.11?}
\expandafter \def \csname FIGLABELHisK\endcsname {3.8?}
\expandafter \def \csname EQLABELequalint\endcsname {3.12?}
\expandafter \def \csname EQLABELtoprove\endcsname {3.13?}
\expandafter \def \csname EQLABELona\endcsname {3.14?}
\expandafter \def \csname FIGLABELhar\endcsname {3.9?}
\expandafter \def \csname FIGLABELeconvex\endcsname {3.10?}
\expandafter \def \csname FIGLABELenested\endcsname {3.11?}
\expandafter \def \csname EQLABELnested\endcsname {3.15?}
\expandafter \def \csname EQLABELeunique\endcsname {3.16?}
\expandafter \def \csname EQLABELsameint\endcsname {3.17?}
\expandafter \def \csname EQLABELequality\endcsname {3.18?}
\expandafter \def \csname FIGLABELsmoothconvex\endcsname {3.12?}
\expandafter \def \csname FIGLABELnotstrict\endcsname {3.13?}
\expandafter \def \csname FIGLABEL2d\endcsname {3.14?}
\expandafter \def \csname EQLABELnorm2\endcsname {4.1?}
\expandafter \def \csname FIGLABELrearrangement\endcsname {4.1?}
\expandafter \def \csname EQLABELbvinequality\endcsname {4.2?}


\hfuzz=5pt
\forwardreferencetrue
\theoremnumber=0
\equationnumber=0
\def\R{{\bf R}}
\def\fracnum#1/#2{\leavevmode
\kern.1em\raise .5ex\hbox{\the\scriptfont0 #1}%
\kern-.1em 
/\kern-.15em\lower.25ex\hbox{\the\scriptfont0 #2}}
\def\mch{{\cal H}_{\partial H}}
\def\mcen{{\cal H}_{\partial E_n}}
\def\mce{{\cal H}_{\partial E}}
\def\mcf{{\cal H}_{\partial F}}
\footline{}
\vskip1.7truein	  
\centerline{\bf Area Minimizing Sets Subject to a Volume Constraint 
in a Convex Set} \vskip1truein
\setbox2=\vtop{
\halign{
#\hfill&\qquad#\hfill\cr
Edward Stredulinsky\footnote{$^{1}$}{Research supported in part 
by a grant from the National Science Foundation}&William P. Ziemer
\footnote{${}^{2}$}{Research supported in part by a grant from
the National Science Foundation}\cr
Mathematics Department&Mathematics Department\cr
University of Wisconsin&Indiana University\cr
Center-Richland&\cr
Richland Center, WI&Bloomington, IN 47405\cr}}
\centerline{\box2}
\vfill\eject
\footline{\centerline{-- \folio --}}
\pageno=1
\sectionnumber=1
{\noindent\bf 1. Introduction.} 
\bigskip
In this paper we consider the problem
of minimizing area subject to a volume constraint in a given convex
set. In precise terms we have the following. Let $ \Omega
\subset\R^{n}$ be a bounded convex set. Thus, $ \abs{\Omega
}<\infty $ where $ \abs{\Omega }$ denotes Lebesgue measure. For a
number $ 0<v<\abs{\Omega }$, let $ E\subset\Omega $ denote a set
with $ \abs{E}=v$ such that 
$$
P(E) \leq P(F)
$$
for all sets $ F\subset\Omega$ with $ \abs{F}=v$, where $
P(E)$
denotes the perimeter of $ E$. The main question we investigate is
whether $ E$ is convex.

It should be emphasized that the perimeter of a competitor $ F$ is
taken relative to $ \R^{n}$, or what is the same, the perimeter is
taken relative to the closure of $ \Omega $ since $ F$ is assumed to
be a subset of $ \Omega $. This problem is considerably
different from minimizing perimeter relative to the interior of $
\Omega $. This was considered in [Gr] where it was shown that a
minimizer is regular and intersects $ \partial \Omega $ orthogonally.

The question of existence of a solution to our problem is resolved
immediately in the context of sets of finite perimeter. Regularity
questions have been considered by other authors. Tamanini [T] has
shown that an area minimizing set $ E$ subject to a volume constraint
has the property that $\partial E\cap\Omega $ is real analytic except
for a closed set whose Hausdorff dimension does not exceed $ n-8$.
Also, under the assumption that $ \partial \Omega \in C^{1}$, it was
shown in [GMT2] that $ \partial E$ is an $ (n-1)$ manifold of class $
C^{1}$ in some neighborhood of each point in  $ \partial
E\cap\partial \Omega $. In $\R^2$, and in $R^n,\,n>2$ 
under an additional condition on $ \Omega$,
we are able to obtain regularity results and ultimately establish that
a minimizer $ E$ is convex. Assuming only that $ \Omega $ is 
bounded and convex, the convexity of $ E$ is an open question in $\R^n,\,n>2$.

The additional condition we impose on $ \Omega $ if $n>2$ is the
following. 
\setbox3=\vbox{\hsize=.8\hsize \noindent 
We assume that a
largest closed ball, $B_{\Omega }$, contained in $ \Omega $ has a
great circle that is a subset of $ \partial \Omega $. 
A great circle of $
B_{\Omega }$ is defined as the intersection of $ \partial B_{\Omega }$
with a hyperplane, $T_{B_{\Omega }}$, passing through the center of $
B_{\Omega }$. The equatorial ``disk'' is defined as 
$D_{B_{\Omega}}=T_{B_{\Omega}}\cap B_{\Omega }$.}
$$
\cases{
\box3&\cr}\eqnlbl{assumption}
$$
Also, assuming initially that $ \partial \Omega \in
C^{2}$ and strictly convex, we invoke a result of [BK] to conclude
that  $\partial E\in C^{1,1}$ at points near $ \partial \Omega $. We
then show, Theorem \thmref{econvex}, that $ E$ is convex. Finally, 
through an approximation procedure, we show that $ E$ is convex with
$ C^{1,1}$ boundary assuming only that $ \Omega $ satisfies a great
circle condition. Clearly, there is no uniqueness if $ v$ is too
small. However, with $ H_{\Omega }$ denoting the union of all largest
balls in $ \Omega $, if $ \abs{H_{\Omega }}\leq v<\abs{\Omega }$, then
$ E$ is unique. In addition for such $v$ we show that perimeter
minimizers $E$ are nested as a function of $v$. In general for nonconvex
$\Omega$ one can expect neither uniqueness or nestedness as indicated
by examples in [GMT1].

The nestedness of perimeter minimizers allows one to rearrange level
sets of functions to create test functions useful in studying
minimizers to certain variational problems. For domains $\Omega$
having certain symmetries it is frequently possible to apply
symmetrization to gain information on minimizers of functionals
such as 
$$
\int_\Omega |\nabla u|^p +\int_\Omega F(u) + \int_0^{|\Omega|} G(u^*,{u^*}')
$$
over appropriate function classes , where $u^*$ is the decreasing
rearrangement of $u$.
However this greatly restricts the collection of domains 
which can be considered.
In Section 4 for the case $p=1$ we construct a rearrangement which 
retains various useful properties of symmetrization while allowing 
a much larger class of domains to be considered, namely those convex domains 
described above.  This rearrangement is useful when 
one has a boundary condition
of the form $u=0$ on $\partial\Omega$ and when it can be established, 
for instance
using truncation, that $u\ge0$ in $\Omega$.
Since this rearrangement produces functions of bounded variation it is more
accurate to replace $\int |\nabla u|$ in the functional above by
the BV norm. The results of Section 3
allow one to deduce certain  regularity properties for minimizers $u$ 
from regularity properties of $u^*$. In addition they establish 
the convexity of 
the sets $\{u>t\}$. Results in [LS]  show that one can not hope
for similar results if $p>1$.
\bigskip
\sectionnumber=2
{\noindent\bf 2. Notation and Preliminaries.} 
\bigskip
The Lebesgue measure of a set $E\subset \R^{n}$ will be denoted by
$\abs{E}$ and $H^{\alpha}(E)$, $\alpha>0,$ will denote its
$\alpha$-dimensional Hausdorff measure. 
If $\Omega\subset \R^{n}$ is an open
set, the class of functions $u\in L^{1}(\Omega)$ whose partial
derivatives in the sense of distributions are measures with finite
total variation in $\Omega$ is denoted by $BV(\Omega)$ and is called
the space of {\it functions of bounded variation in $\Omega.$} The
space $BV(\Omega)$ is endowed with the norm 
$$
\norm{u}_{BV(\Omega)}=\norm{u}_{1;\Omega}+\norm{\nabla u}(\Omega)\eqnlbl{np.1}
$$
where $\norm{u}_{1;\Omega}$ denotes the $L^{1}$-norm of $u$ on
$\Omega$ and 
where $\norm{\nabla u}$ is the total variation of the vector-valued
measure $\nabla u$.

A 
Borel set $E\subset \R^{n}$ is said to have {\it finite perimeter in
$\Omega$} provided the characteristic function of $E$, $\upchi_{E}$, is a
function of bounded variation in $\Omega$. Thus, the partial
derivatives of $\upchi_{E}$ are Radon measures on $\Omega$ and the
perimeter of $E$ in $\Omega$ is defined as
$$
P(E,\Omega)=\norm{\nabla\upchi_{E}}(\Omega).\eqnlbl{np1}
$$
A set $E$ is said to be of
{\it locally finite perimeter} if $P(E,\Omega)<\infty$ for every
bounded open set $\Omega\subset \R^{n}$. 

The definition implies that sets of finite perimeter are defined only
up to sets of measure 0. In other words, each set determines an
equivalence class of sets of finite perimeter.  In order to avoid
this ambiguity, whenever a set $E$ of finite perimeter is considered
we shall always employ the measure theoretic closure as the set to
represent $E$. Thus, with this convention, we have
$$
x\in E\;\hbox{if and only if}\;\limsup_{r\to0}\frac{\abs{E\cap
B(x,r)}}{\abs{B(x,r)}}>0. \eqnlbl{rep}
$$

One of the fundamental
results of the theory of sets of finite perimeter is that they
possess a measure-theoretic exterior normal which is suitably general
to ensure the validity of the Gauss-Green theorem. A unit vector
$\nu$ is defined 
as the exterior normal to $E$ at $x$ provided
$$
\lim_{r\to0}r^{-n}\abs{B(x,r)\cap\{y:(y-x)\cdot\nu<0,y\notin E\}}=0
$$
and
$$
\lim_{r\to0}r^{-n}\abs{B(x,r)\cap\{y:(y-x)\cdot\nu>0,y\in E\}}=0,\eqnlbl{np2}
$$
where $B(x,r)$ denotes the open ball of radius $r$ centered at $x$.
The measure-theoretic normal of $E$ at $x$ will be denoted by
$\nu(x,E)$ and we define
$$
\partial^{*}E=\{x:\nu(x,E)\;{\rm exists}\}.\eqnlbl{np3}
$$
Clearly, $\partial^{*}E\subset\partial E$, where $\partial E$ denotes
the topological boundary of $E$.

A set $E$ of finite perimeter is said to be {\it area minimizing in
an open set} $\Omega$ if
$$
\norm{\nabla \upchi_{E}}(\Omega)\leq\norm{\nabla
\upchi_{F}}(\Omega)\eqnlbl{areamin} 
$$
for every set $F$ with $F\Delta E\subset\subset \Omega$. Here
$F\Delta E$ denotes the symmetric difference.

The regularity of $\partial E$ will
play a crucial role in our development. Suppose $\partial E$ is
area minimizing in $ U$
and for convenience of notation, suppose $0\in U\cap\partial E$. For each
$r>0$, let $E_{r}=\R^{n}\cap\{x:rx\in E\}$. It is known (cf.
[S,\S 35], [MM,\S 2.6]) 
that for each sequence $\{r_{i}\}\to0$ there exists a subsequence
(denoted by the full sequence) such that $\upchi_{ E_{r_{i}}}$ converges
in $L^{1}_{\rm loc}(\R^{n})$ to $\upchi_{C}$, where $C$ is a set of locally
finite perimeter. In fact, $\partial C$ is area minimizing and is called the
tangent cone to $E$ at 0. Although it is not immediate, $C$ is a cone
and therefore the union
of half-lines issuing from 0. It follows from [S, \S 37.6] that if
$\overline{C}$ is contained in $\overline{H}$ where $H$ is any
half-space in $\R^{n}$ with $0\in\partial H$, then $\partial E$ is
regular at 0. That is, there exists $r>0$ such that 
$$
B(0,r)\cap\partial E\;\,\hbox{is a real analytic hypersurface.}\eqnlbl{np10}
$$
Furthermore, $\partial E$ is regular at all points of $\partial
^{*}E$ and 
$$
H^{\alpha}((\partial E - \partial ^{*}E)\cap U)=0
\quad\hbox{\rm for all}\;\alpha>n-8,\eqnlbl{np11}
$$
cf. [Gi, Theorem 11.8].

The notion of {\it excess} plays a critical role in the theory of
minimal boundaries. It measures how far a set $E$ is from being area
minimizing in a ball. Formally, it is defined by
$$
\psi(x,r)=\norm{\nabla \upchi_{E}}(B(x,r))-\inf\{\norm{\nabla
\upchi_{F}}(B(x,r)):F\Delta E\subset\subset B(x,r)\}. 
$$
Thus, $\psi\equiv0$ when $E$ is area minimizing. If $E$ is an
arbitrary set of finite perimeter and $\psi(x,r)\leq
Cr^{n-1+2\alpha}$ for some $x\in\partial E$ and all $0<r<R$ with given
constants $C, R$ and $0<\alpha<1$, then it follows from a result of
Tamanini [T] that there is an area minimizing tangent cone to
$\partial E$ at $x$. 
\bigskip
{\bf \thmlbl{lipdiff} Definition.} Let $M$ denote a $k$-dimensional
$C^{1}$ submanifold of $\R^{n},\>0<k<n,$ and let $f\colon M\to\R$ be
an arbitrary function. We will say that $f$ is differentiable at $x_{0}\in M$
if $f$ is the restriction to $M$ of a function $\bar f\colon U\to\R$
where is $U\subset\R^{n}$ is some open set containing $x_{0}$ and
where $\bar f$ is differentiable at $x_{0}$.
\bigskip
{\bf \thmlbl{diffae} Lemma.} {\it Let $M$ be an $n-1$-dimensional
$C^{1}$ submanifold of $\R^{n}$ and let $f\colon M\to\R$ be a
Lipschitz function. Then $f$ is differentiable at $H^{n-1}$ almost
all points of $M$.}
\medskip
{\bf Proof.} The manifold $M$ near any of its points $x_{0}$ can be
represented as the graph of a function defined on some open $n-1$-ball
$B'\subset\R^{n-1}$. Thus, there is an open $n$-cylinder $C$ of the
form $C=B'\times(a,b)$ such that $C-M$ consists of two nonempty 
connected, open
sets and that each projection of $M\cap C$ onto the top and
bottom of $C$ is a homeomorphism. Let points $x\in C$ be denoted by
$x=(x',y)$ where $x'\in B'$ and $y\in(a,b)$ and define $\bar f\colon
C\to\R$ by $\bar f(x',y)=f(x',y_{M})$ where $(x',y_{M})$ is that
unique point on $M\cap C$ that is the projection of $(x',y)$. It is easy to
verify that $\bar f\colon C\to\R$ is Lipschitz and therefore, by
Rademacher's theorem, that $\bar f$ is differentiable at (Lebesgue)
almost all points of $C$. Let $N$ denote those points at which $\bar
f$ is not differentiable. Clearly, if  $\bar f$ is
differentiable at a point $(x',y_{1})$ then it is
differentiable at any other point of the form $(x',y_{2})$.
Now define $d\colon C\to\R$ by
$d(x',y)=\abs{y-y_{M}}$. Note that $d$ is Lipschitz and that
$d^{-1}(t)$ consists of two 
copies of $M\cap C$, one is a vertical distance of $t$ units above
$M\cap C$ and the other is a vertical distance of $t$ units below $M\cap
C$. Now employ the co-area formula to obtain
$$
0=\int_{N}\abs{\nabla d}\;dx=\int_{a}^{b}H^{n-1}[d^{-1}(t)\cap N]\;dt.
$$
Thus, for almost every $t\in(a,b),\>\bar f$ is differentiable at
$H^{n-1}$ almost all points of $d^{-1}(t)$. Consequently, $\bar f$ is
differentiable at the corresponding points of $d^{-1}(0)=M\cap C$;
that is, $\bar f$ is differentiable at $H^{n-1}$ almost all points of
$M\cap C$, as required.\qed
\bigskip
In view of the preceding Lemma, we can define the directional
derivative of $f$ relative to $M$ at $H^{n-1}$-almost all $x\in M$ in
the usual manner. Given a vector $\tau$ in the tangent space to
$M$ at $x$, let $\gamma\colon(-1,1)\to M$ be any $C^{1}$ curve with
$\gamma(0)=x$ and $\gamma'(0)=\tau$. Define 
$$
D_{\tau}f(x)=(\bar f\circ\gamma)'(0)
$$
where it is understood that $\bar f$ is differentiable at $x$.
Observe that this definition is independent of the extension $\bar
f$. 

If we are given a Lipschitz vector field $X\colon M\to\R^{n}$, by using
usual methods, it now becomes clear how to define the divergence of
$X$ relative to $M$, denoted by $\div_{M}X$.

If the closure $\overline{M}$ of $M$ is a $C^{1}$ manifold with boundary
$\partial M=\overline{M}-M$ and if $X\colon\R^{n}\to\R^{n}$ is a $C^{1}$
vector field with the property that for each $x\in M,\;X(x)$ is an
element of the tangent space to $M$ at $x$, then the classical
divergence theorem states 
$$
\int_{M}\div_{M}X\;dH^{n-1}=\int_{\partial M}X\cdot\eta
\;dH^{n-2}\eqnlbl{div} 
$$
where $\eta$ is the outward pointing unit co-normal of $\partial M$.
That is, $\abs{\eta}=1$, $\eta$ is normal to $\partial M$, and
tangent to $M$.
\bigskip
{\bf \thmlbl{meancur} Definition.} Let $M$ be an
oriented $n-1$-dimensional
submanifold of $\R^{n}$ of class $C^{1,1}$; that is, $M$ is of class
$C^{1}$ and its unit normal $\nu$ is Lipschitz. From Lemma
\thmref{diffae}, we have that the components of $\nu$ are
differentiable at $H^{n-1}$ almost all points of $M$. Thus,
$\div_{M}\nu$ is defined $H^{n-1}$ almost everywhere on $M$. At such
points, we
define the {\it mean curvature} of $M$ at $x$ as 
$$
{\cal H}_{M}(x)=\div_{M}\nu(x)
$$
If $X\colon\R^{n}\to\R^{n}$ is a $C^{1}$ vector field, consider its
decomposition into its tangent and normal parts relative to $M$,
$$
X=X^{\top}+X^{\bot}
$$
where
$$
X^{\bot}=(X\cdot\nu)\nu.
$$
Then, at $H^{n-1}$ almost all points in $M$, it follows that
$$
\div_{M}X^{\bot}=(X\cdot\nu)\div_{M}\> \nu.
$$
Hence,
$$
\div_{M}X^{\bot}={\cal H}_{M}X\cdot\nu.
$$
On the other hand, from \eqnref{div} we have
$$
\int_{M}\div_{M}X^{\top}\;dH^{n-1}=\int_{\partial
M}X\cdot\eta\;dH^{n-2}.
$$
Since $\div_{M}X=\div_{M}X^{\top}+\div_{M}X^{\bot}$, we obtain
$$
\int_{M}\div_{M}X\;dH^{n-1}=\int_{M}{\cal
H}_{M}X\cdot\nu\;dH^{n-1}+\int_{\partial M}X\cdot\eta\;dH^{n-2}.\eqnlbl{divmc}
$$
\sectionnumber=3
\theoremnumber=0
\equationnumber=0
{\noindent\bf 3. Main Results}
\bigskip
In this section we consider the following situation. 
\setbox2=\vbox{\hsize=.8\hsize
\noindent Let $\Omega$ be a
bounded,  convex domain in $\R^{n},\>n\geq2$.
 Let $E\subset\overline\Omega$ denote a set which
minimizes perimeter in the closure of $\Omega$ subject to a volume
constraint $\abs{E}=v<\abs{\Omega}$. Thus 
$$
P(E,\R^{n})\leq P(F,\R^{n})
$$
for all sets $F\subset\overline\Omega$ with $\abs{F}=v$.}
$$
\cases{
\box2&\cr}\eqnlbl{problem}
$$
\smallskip
We will first establish boundary regularity and curvature properties
for such perimeter minimizers under the assumption that $\Omega$
is {\it strictly convex and that} $\partial \Omega\in
C^{2}$. Convexity, nestedness and uniqueness results will then be established
under the further assumption that
$$
n=2 \quad\hbox{or}\quad \Omega \hbox{  satisfies a great circle
condition. }
$$
The assumption of strict convexity and $C^2$ regularity will
then be dispensed with in part through an approximation argument.

Associated with \eqnref{problem} is some further notation. We let $H$
denote the convex hull of a minimizer $E$ of \eqnref{problem}, and we
denote by $H^{+}$ that part of $H$ that lies ``above'' the equatorial 
disk $D_{B_{\Omega }}$ of $B_{\Omega}$ as defined in
\eqnref{assumption}. Since
$P$
divides
$H$
into
two parts, we arbitrarily
call one of them the part that lies ``above'' $P$. 

Next, we recall some facts concerning area minimizing sets with a
volume
constraint. The main result of [GMT1] is that if $E$ is area
minimizing with a volume constraint, then 
$$
\psi(x,r)\leq Cr^{n}\eqnlbl{excessrn}
$$
for each $x\in\partial E$ and for all sufficiently small $r>0$.
Consequently, it follows from work of Tamanini [T] that an area
minimizing set $E$ with a 
volume constraint possesses an area minimizing tangent cone at each
point of $(\partial E)\cap\Omega$.
From this it follows that $(\partial E)\cap\Omega$ enjoys
the same regularity properties as an area minimizing set; that is,
$(\partial E)\cap\Omega$ is real analytic except for a closed
singular set $S$ whose Hausdorff dimension does not exceed $n-8$.
Furthermore, it was established in [GMT2, Theorem 3] that $\partial
E$ is an $(n-1)$ manifold of class $C^{1}$ in some neighborhood of
each point $x\in\partial E\cap\partial \Omega$. 
 
The object of this section is to prove that $E$ is convex and we
begin by proving $C^{1,1,}$ regularity of $\partial E$ near
$\partial \Omega$. For this we will need the following result of
Br\'ezis and Kinderlehrer, [BK]. 
\bigskip
{\bf \thmlbl{bk} Theorem.} {\it Let $a\colon\R^{n-1}\to\R^{n-1}$ be a
$C^{2}$ vector field satisfying the condition that for each compact
$C\subset\R^{n-1}$, there exists a constant $\nu=\nu(C)>0$ such that 
$$
(a(p)-a(q))\cdot(p-q)\geq\nu\abs{p-q}^{2}
$$
for all $p,q\in C$. Let $U\subset\R^{n-1}$ be an open connected set
and let $\beta\in C^{2}(U)$ satisfy
$\beta\leq0$ on $\partial U$. Let $f\in C^{1}(U)$. With
$\bf K$=${\bf K}_{\beta}$ denoting the convex set of Lipschitz
functions $v$ satisfying $v\geq\beta$ in $U$ and $v=0$ on $\partial
U$, let $u\in{\bf K}$ be a solution of 
$$
\int_{U}a(\nabla u)\cdot\nabla (v-u)\;dx\geq\int_{U}f(v-u)\;dx
$$
for all $v\in {\bf K}$. Then $u\in C^{1,1}(V)$ on any domain $V$ with
$\overline{V}\subset U$.}
\bigskip
We now apply this result to obtain $C^{1,1}$ regularity of the
boundary of a minimizer $E$ of the variational  problem
\eqnref{problem} near $\partial \Omega$ . Since $\partial E$ is an
$(n-1)$ manifold of class $C^{1}$ in some neighborhood of each point
$x\in\partial E\cap\partial \Omega$, it follows that near such a point
$x$, we may represent both $\partial E$ and $\partial \Omega$ as
graphs of functions $u$ and $\beta$, respectively, defined on an open
set $U'\in\R^{n-1}$ containing $x'$ where $x=(x',y''),\;y''\in\R$.
We
will assume $u$ and $\beta$ chosen in such a way that
$u\geq\beta,\;u=0$ on $\partial U'$ and $\beta\leq0$ on $\partial U'$.
Using the  convexity of $\Omega$, this can be accomplished by
considering a hyperplane $P_{0}$ passing through $E$ and parallel to
the tangent plane to $\partial E$ at  $x$. By taking $P_{0}$
sufficiently close to the tangent plane, $U'$ can be defined as
$P_{0}\cap E$. Now select $v\in {\bf K}$ and for $0<\varepsilon<1$,
define $u_{\varepsilon}$ on $U'$ as
$u_{\varepsilon}=u+\varepsilon(v-u)$. We will assume $\varepsilon$
chosen small enough so that the graph of $u_{\varepsilon}$ remains in
$\overline{\Omega}$. Note that $u_{\varepsilon}\in {\bf K}$. Select a
point $z\in(\partial E)\cap\Omega$ at which $\partial E$ is regular.
Thus, $\partial E$ is real analytic near $z$ and its mean curvature is
a constant $K$ there. In a neighborhood of $z$, we can represent
$\partial E$ as the graph of a function $w$ defined on
some open set $V'\subset\R^{n-1}$ containing 
$z'$ where $z=(z',z'')$. The neighborhoods about $x$ and $z$ where
$\partial E$ is represented as a graph are taken to be disjoint.
Let $\varphi\in C^{\infty}_{0}(V')$ denote a
function with the property that
$$
\int_{V'}\varphi\;dH^{n-1}=\int_{U'}(v-u)\;dH^{n-1},\eqnlbl{samevol}
$$
and define $w_{\varepsilon}=w-\varepsilon\varphi$. The graphs of the functions
$u_{\varepsilon}$ and $w_{\varepsilon}$ produce a perturbation of the
set $E$, say $E_{\varepsilon}$. Because of \eqnref{samevol}, we have
that $\abs{E}=\abs{E_{\varepsilon}}$. With 
$$
F(\varepsilon)=\int_{U'}\sqrt{1+\abs{\nabla u_{\varepsilon}}^{2}}+
              \int_{V'}\sqrt{1+\abs{\nabla w_{\varepsilon}}^{2}},
$$
the minimizing property of $\partial E$ implies that $F(0)\leq
F(\varepsilon)$ for all small $\varepsilon$ and therefore that
$F'(0)\geq0$. Thus,
$$
\int_{U'}\frac{\nabla u}{\sqrt{1+\abs{\nabla u}^{2}}}\cdot\nabla
(v-u) -\int_{V'}\frac{\nabla w}{\sqrt{1+\abs{\nabla w}^{2}}}\cdot\nabla
\varphi\geq0.
$$
Since $w$ has constant mean curvature $K$, we obtain
$$
\int_{V'}\frac{\nabla w}{\sqrt{1+\abs{\nabla w}^{2}}}\cdot\nabla
\varphi=-\int_{V'}K\varphi=-K\int_{V'}\varphi=-K\int_{U'}(v-u),
$$
and therefore
$$
\int_{U'}\frac{\nabla u}{\sqrt{1+\abs{\nabla u}^{2}}}\cdot\nabla
(v-u) \geq -K\int_{U'}(v-u).\eqnlbl{mceleqk}
$$
\qed

If $\eta\in C^{\infty}_{0}(U')$ denotes an arbitrary nonnegative
test function, then with $v-u=\eta$, \eqnref{mceleqk} states that
$u$ is a weak solution of $\mce\leq K$. This combined with the $C^{1,1}$-
regularity of $u$ implies that $\mce\le K$ pointwise almost everywhere
in a neighborhood of $\partial \Omega$. Since $\mce=K$ in $\partial
E\cap(\Omega\setminus S)$ 
with $H^{n-1}(S)=0$ we have the following result.
\bigskip
{\bf\thmlbl{beregular} Theorem.} {\it Assume that $\Omega$ is bounded,
convex and has a $C^2$ boundary. If $E$ is a minimizer of
{\rm \eqnref{problem}}, then
$\partial E\in C^{1,1}$ in some neighborhood 
of $\partial \Omega$ and $\mce\le K$ $H^{n-1}$-almost everywhere on
$\partial E$.} 
\bigskip
We now will exploit Theorem \thmref{beregular} to establish both regularity
and a mean curvature estimate for the boundary of the convex hull of $E$.

\medskip
{\bf\thmlbl{bhregular} Theorem.} {\it Assume that 
$\Omega$ is bounded, strictly
convex and has a $C^2$ boundary. If $E$ is a minimizer of
{\rm \eqnref{problem}} with convex hull $H$ then
$\partial H\in C^{1,1}$  and $\mch\le K \;H^{n-1}$-almost everywhere
on $\partial H$.} 
\medskip
{\bf Proof.}
Note that the singular set $S$ in $\partial E$ is a 
closed subset of $\Omega$ and thus separated from $\partial\Omega$, in fact it
is contained in the interior of H, for if 
$x\in\partial E\cap\partial H\cap\Omega$,
then the tangent cone to $\partial E$ at $x$ must be a hyperplane
because $E\subset H$ and $H$ is convex.  Consequently $\partial E$ is
regular at $x$. Let $N$ be an open 
neighborhood of $S$ with compact closure in the interior of $H$. Thus by 
Theorem \thmref{beregular} and the analyticity of $\partial E$ in
$\Omega\setminus S$ 
we see that
$\partial E$ is $C^{1,1}$ at points in $G:=\partial E\setminus
N$.
Therefore for some $C$ we have
$$
|\nu(x)-\nu(z)|\le C|x-z|\qquad\qquad x, z\in G\eqnlbl{lip}
$$
where $\nu(x)$ is the outward unit normal to $\partial E$ at $x$.
Also since $\partial E$ is $C^1$ at points in $G$ there exists
an $\varepsilon$ such that for all $x\in G$ and $z\in\partial E\cap
B(x,\varepsilon)$ we have
$$
|\nu(x)\cdot(x-z)|\le\fracnum{1}/{2}|x-z|.\eqnlbl{c1}
$$

Choose $x\in\partial E\cap\partial H\subset G$ and let $ 0<\alpha
<1/2$. Then define 
$$
d=\alpha \min\{\varepsilon, \hbox{ dist}(\partial H, N), {1\over
2C},\hbox{\rm diam}\>E\}.
$$ 
Let $ y=x-d\nu (x)$ and observe that
$ y$ is in the interior of $ E$ since $\partial E$ cannot intersect the
line segment $\overline{xy}$ at a point $z\not=x$ due to \eqnref{c1} . Let
$r=\hbox{ dist}(y,\partial E)$ and note that $0<r\le d$. Now choose
any 
$z\in\partial E$ such that $|y-z|=r$. Note that $z\in G$, for
otherwise we would have $ z\in N$ and since $
\abs{x-z}\leq\abs{x-y}+\abs{y-z}$, it would follow that 
$$
2d\geq\abs{x-z}\geq\hbox{\rm dist}(\partial H,N)\geq\frac{d}{\alpha
}>2d,
$$
a contradiction. 
Then, $|x-z|\le|x-y|+|y-z|\le2 d<\varepsilon$ and both \eqnref{lip}
and \eqnref{c1} hold. Thus, since $x=y+d\nu(x)$ and $z=y+r\nu(z)$,
we have
$|d-r|\le|\nu(x)\cdot(x-z)|$ and 
$$
|x-z|=|(d-r)\nu(x)+r(\nu(x)-\nu(z))|\le(\fracnum{1}/{2}+C
r)|x-z|\le\fracnum{3}/{4}|x-z|,
$$
(since $r\le d\le \alpha/(2c)\le 1/(4c)$)
which implies that $x=z$ and therefore $r=d$. This implies that
for every $x\in \partial E\cap\partial H$ there exists 
a ball $B_x\subset E$ of radius $d$ containing $x$.

Given any $p\in \partial H$
we claim that $p$ is a convex combination of points $\{x_i\}$ in
$\partial E\cap\partial H$. To see this note that if $C$ is a
convex set with $E\subset C$ then $\overline 
E\subset C$ since if $x\in \overline E$
then either $x\in C$ or $x\in\partial C$; in the later case $x$ lies in a 
support plane of $C$ so if $x\in \Omega$, regularity theory implies
that $x\in E\subset C$, and if $x\in \partial\Omega$ then $x$ is not in the
singular set $S$ of $E$ (since $S$ is a compact subset of $\Omega$)
so again $x\in E\subset C$. Consequently from the definition of convex
hull $H$ of $E$ as the intersection of all convex sets containing
$E$, we see that $\overline E\subset H$. Moreover $H$ is the convex hull of
$\overline E$ 
from which we conclude by a well known result that $H$ is closed since
$\overline E$ is a compact subset of $\R^n$. Note that the set of finite
convex combinations of points from $E$ is convex, contains $E$, and
is contained in any convex set which contains $E$ and so equals $H$.
Thus if $p\in\partial H$ we have  $p\in H$, since $H$ is closed,
and consequently $p=\sum_{i=1}^k \lambda_i x_i$ for $x_i\in E$ and 
$\sum_{i=1}^k \lambda_i=1,\,
 \lambda_i\ge0, i=1\dots k$. If we take $k$ to be as small as possible
then either $k=1$ and $p\in E$ and the claim is trivially true, or
$p$ lies in the $k$ dimensional interior of the convex hull $M$ of 
$\{x_i\}$ in which case no $x_i$ can lie in the interior of $H$ since then 
the same would be true of $p$. Consequently $x_i\in \partial
E\cap\partial H,\,i=1\dots k$, as claimed.

Taking  the convex hull of $\cup_{i=1}^k B_{x_i}$ we see that there
exists a ball $B_p\subset H$ 
of radius $d$ containing $p$, i.e. $H$ satisfies a uniform
interior sphere 
condition.  We claim that this implies $\partial H$ is $C^{1,1}$.
To see this, consider the problem of prescribing unit vectors $\nu_1,
\nu_2\in\R^n$, and finding a convex set $\tilde H$, 
satisfying the interior sphere condition noted above, and points $x,
y\in\partial\tilde H$ 
with $\nu(x)=\nu_1, \nu(y)=\nu_2$, such that $|x-y|$ is minimized.
It is clear that $x, y$ must lie in a two 
dimensional plane orthogonal to the intersection of  two
hyperplanes having $\nu_1, \nu_2$ as normals, i.e. one need only
consider the two dimensional 
case where it is easy to see that one must have $B_x=B_y$.
Taking the center of this ball to be the origin then
$\nu(x)=x/d, \nu(y)=y/d$ and  
we trivially have
$$
|\nu(x)-\nu(y)|\le{1\over d}|x-y|.
$$
Since this is the case when $|x-y|$ is smallest for fixed $\nu(x), \nu(y)$
we have established that $\nu(x)$ is Lipschitz in general.

We now  prove that $\mch\leq K\;H^{n-1}$-almost everywhere in
$\partial H$. 
Note that $\mch=\mce$ $H^{n-1}$-almost everywhere
on $\partial E\cap\partial H$ by Theorem \thmref{beregular}. Thus we
need only consider 
points $p\in\partial H\setminus\partial E$.
In fact since $\partial H$ is $C^{1,1}$ we need only consider
$p\in\partial H\setminus\partial E$ at which $\partial H$ is
classically twice differentiable. As above, any such $p$ lies 
in the $k$ dimensional interior of the convex hull $M$ of certain points 
$p_i\in \partial E,\, i=1,\dots, k$. Note that $k\ne 1$ due to
$p\notin\partial E$.
Choose a coordinate system such that points in $\R^n$ are represented
as $(x,y,z),\,x\in\R^k,\, y\in\R^{n-k-1},\, z\in\R$, with $z=0$ the
tangent plane to $\partial H$  at $p$, $p_i=(x_i,0,0),\, i=1,\dots,
k$, and $z\ge 0$ in $H$. We will construct an
analytic function $g$ whose graph does not lie below $\partial H$,
contains $M$, 
and has mean curvature bounded above by $K+\varepsilon$ (for any
$\varepsilon>0)$ in a small neighborhood of $p$. This will lead to
the conclusion that $\mch\le K$ at $p$.

Let $\partial E$ be represented as $z=f(x,y)$ for $f$ defined in a
neighborhood in $\R^{k}\times\R^{n-k-1}$ of $\cup(x_{i},0)$. 
Thus
$$
(x_i,y, f(x_i,y))\in \partial E\subset H
$$
for small $|y|$, and consequently
$$
\sum_{i=1}^k \lambda_i(x_i,y,f(x_i,y)) \in H \qquad \hbox{ if }\quad
\sum_{i=1}^k \lambda_i=1,\, \lambda_i\ge0 \eqnlbl{suminh}
$$
for small $|y|$.
For any given $x$ in $N$, where $N$ is the convex hull of the points
$x_i,\,i=1,\dots, k$, 
let $\lambda=\lambda(x)=(\lambda_1(x),\dots, \lambda_k(x))$
be the unique vector such that 
$$
x=\sum_{i=1}^{k} \lambda_i(x) x_i,\,\quad
\sum_{i=1}^k \lambda_i(x)=1,\,\quad \lambda_i(x)\ge0.
$$
Thus if we define 
$$
g(x,y)=\sum_{i=1}^k \lambda_i(x) f(x_i,y)
$$
we see from \eqnref{suminh} for $x\in N$ and small $|y|$ that 
$$
(x,y,g(x,y))\in H,
$$
and so the surface $z=g(x,y)$ does not lie below $\partial H$
at such $(x,\,y)$.

Note that $M\cap\partial \Omega =\nullset$, for otherwise the
plane $z=0$, which contains $M$, would be a tangent plane to $\partial
\Omega $, thus contradicting the strict convexity of 
$\partial \Omega$.
Also $M$ does not intersect the 
singular set of $\partial E$ since $M\subset \partial H$. Thus 
$\partial E$ is analytic at each $p_i$ and therefore both $f(x_i,y)$
and $g(x,y)$ are smooth for small $|y|$. Furthermore, 
$$
0\le \Delta_y f(x_i,0) \le \Delta f(x_i,0) \le   K
$$
since $\nabla f(x_i,0)=0$, $\mce$ equals $  \Delta f$  at
points where the gradient is zero, and the second derivatives of $f$  are
nonnegative at $(x_i,0)$ due to the fact that $f\ge 0$,
$f(x_i,0)=0$ for all $i$.  Hence, for any $\varepsilon>0$,
$\Delta_y f(x_i,y)\le  (K+\varepsilon)$ for small enough $|y|$
so $\Delta_y g(x,y)\le  (K+\varepsilon)$ as well. However
$\Delta_x g=0$ and so $\Delta g \le  (K+\varepsilon)$
for small $|y|$. Recall that $\partial H$ is trapped between
$\{z=0\}$ and the graph of $g$ over a region which contains $p$
in its interior. Since $g(p)=0$ and $\partial H$ is twice 
differentiable at $p$ we conclude that $\mch(p)\le K$  as required.
\qed 

\bigskip

{\bf\thmlbl{ballinmin} Theorem.} {\it Assume
that $\Omega$ is bounded, strictly convex and satisfies
a great circle condition. If $E$ is a minimizer of
{\rm \eqnref{problem}} with $|B_\Omega|\le|E|$  then 
$$
B_\Omega\subset E
$$
where $B_\Omega$ is the largest ball in $\Omega$.} 
\bigskip
{\bf Proof.} If $|E|=|B_\Omega|$ then clearly $E$ must be a ball.
Since there is only one largest ball in $\Omega$ due to strict
convexity, we have $E=B_\Omega$.  Otherwise $|B_\Omega|<|E|$. In this
case translate the upper and lower hemispheres of $B_\Omega$ by a
distance $d$ in opposite directions orthogonal to 
$T_{B_{\Omega }}$ until $H$, the convex
hull of the two translated hemispheres, intersects $E$ in a set of
measure $|B_\Omega|$ i.e.
$$
|H\cap E|=|B_\Omega|.\eqnlbl{intersectionisball}
$$ 
This is possible because of the great circle conition and because
$\Omega $ is bounded and convex. Now translate the
hemispheres back to their original positions 
while rigidly carrying along the  parts of $E$ lying in the
exterior of $H$.
Let $\tilde E$ be the union of the translated parts of $E$ with $B_\Omega$.
Note that
$$
|\tilde E|=|E|\quad\hbox{\rm and therefore}
\quad P(\tilde E)\geq P(E).\eqnlbl{esamemeasure}
$$ 
Using a standard inequality, cf. [MM], we have
$$
P(E)+P(H)\ge P(E \cap H)+P(E \cup H)
$$
where $P(S)$ denotes $P(S,\,\R^n)$. For brevity, write 
$D=D_{B_{\Omega }}$. Observe that
$$
P(H)=2dH^{n-2}(\partial D)+P(B_\Omega), \quad P(E \cup H)=P(\tilde
E)+2dH^{n-2}(\partial D)
$$
and thus
$$
P(E)+P(B_\Omega)\ge P(E \cap H)+P(\tilde E).
$$
In view of \eqnref{esamemeasure} it follows that 
$P(E\cap H)\leq P(B_{\Omega })$. But then the isoperimetric
inequality and \eqnref{intersectionisball} imply that $E\cap H$
is a ball. However $\Omega$ contains only one largest
ball and so we must have $E\cap H=B_\Omega$, i.e. $B_\Omega\subset
E$.\qed
\bigskip
Suppose $M$ is an oriented $(n-1)$-dimensional $C^{1}$ submanifold of
$\R^{n}$ and $f\colon M\to\R^{n-1}$ a $C^{1}$ mapping. Let $Jf(x)$
denote the Jacobian of $f$ at $x$ and note that the sign of the
Jacobian depends on the orientation of $M$. We recall the following
result, cf. [Fe, Theorem 3.2.20]: For any $H^{n-1}$-measurable set
$E\subset M$ and any $H^{n-1}$-measurable function $\varphi$,
$$
\int_{E}\varphi[f(x)]\abs{Jf(x)}\;dH^{n-1}(x)=\int
\varphi(y)N(f,E,y)\;dy\eqnlbl{intmul} 
$$
where $N(f,E,y)$ denotes the number (possibly infinite) of points in
$f^{-1}(y)\cap E$. Here equality is understood in the sense that if
one side is finite, then so is the other. In our application
\eqnref{sumdeg} below, we will know the left side is finite,
therefore ensuring that $N(f,E,y)$ is finite for almost all $y$.
\bigskip
{\bf \thmlbl{density} Lemma.} {\it There is a constant $C=C(n)$
such that for each $x\in(\partial E)\cap\Omega$ we have
$$
\frac{H^{n-1}((\partial E)\cap B(x,r))}{r^{n-1}}\leq C
$$ 
for almost all sufficiently small $r>0$.}
\bigskip
{\bf Proof.} It follows from \eqnref{excessrn} that we may as well
assume $\partial E$ is area minimizing. In this case the result
follows immediately from the fact that 
$$
\frac{H^{n-1}((\partial E)\cap B(x,r))}{r^{n-1}}
$$
is nondecreasing in $r$, for $r>0$ sufficiently small, cf. [Fe,
Theorem 3.4.3].\qed
\bigskip\noindent
{\bf \thmlbl{forms} Lemma.} {\it For every $\varepsilon>0$ and any
open set $V\subset\R^{n}$ containing the singular set $S$ of
$\partial E$, there exists an open set $W$ and a Lipschitz
function $f$ such that 
$$
\displaylines{
S\subset W\subset\{f=1\}\cr
\hbox{\rm spt}\;f\subset V\cr
\int_{\partial E}\abs{\nabla f}\;dH^{n-1}\leq \varepsilon.\cr}
$$}
\bigskip
{\bf Proof.} Let $V$ be any open set containing $S$ and let
$\delta=1/2(\hbox{\rm dist}\;S,\R^{n}-V)$. Since $H^{n-7}(S)=0$
and $S$ is compact, there is a finite collection of
open balls $\{B(x_{i},r_{i})\}_{i=1}^{m}$ such that $2r_{i}<\delta,
B(x_{i},r_{i})\cap S\not=\nullset,
S\subset\cup_{i=1}^{m}B(x_{i},r_{i})$ and 
$$
\sum_{i=1}^{m}r_{i}^{n-7}< \frac{\varepsilon}{C},
$$
$C$ as in Lemma \thmref{density}. We will assume that
each ball $B(x_{i},r_{i})$ has been chosen so
that $r_{i}<1$ and that $2r_{i}$ satisfies Lemma \thmref{density}.
Let $W$ denote the 
union of these balls and define 
$f_{i}$ by 
$$
f_{i}(x)=\cases{
1&if $\abs{x-x_{i}}\leq r_{i}$\cr
2-\frac{\abs{x-x_{i}}}{r_{i}}&if
$r_{i}\leq\abs{x-x_{i}}\leq2r_{i}$\cr 
0&if $2r_{i}\leq\abs{x-x_{i}}$.\cr}
$$
In view of Lemma \thmref{density}, it follows that 
$$
\int_{B(x_{i},r_{i})\cap\partial E}\abs{\nabla f_{i}}\;dH^{n-1}\leq
Cr_{i}^{n-2}<Cr_{i}^{n-7}.
$$
Now let $f:=\max_{1\leq i\leq m}f_{i}$. Then $f$ is Lipschitz,
$W\subset\{f=1\},\;\hbox{\rm spt}\;f\subset V$ and 
$$
\eqalign{
\int_{\partial E}\abs{\nabla
f}\;dH^{n-1}&\leq\sum_{i=1}^{m}\int_{B(x_{i},r_{i})\cap\partial
E}\abs{\nabla f_{i}}\;dH^{n-1}\cr
&<C\sum_{i=1}^{m}r_{i}^{n-7}<\varepsilon.\cr}  
$$
\qed
\bigskip
\noindent
\bigskip
{\bf \thmlbl{approxid} Lemma.} {\it Let $T$ denote the $(n-1)$-rectifiable
current determined by $(\partial E)^{+}$, the part
of $\partial E$ that lies above the equatorial disk 
$D:=D_{B_{\Omega}}$ of $B_{\Omega}$. Then $\partial T$ is the
$n-2$-sphere given by $\partial T=\partial D$.}
\bigskip
{\bf Proof.} Clearly, the support of $\partial T$ contains the
$n-2$-sphere, but we must rule out the possibility of it containing
points of $S$ as well. For this purpose, choose $x\in S$
and let $\varphi$ be any smooth differential form supported in some
neighborhood of $x$ that does not meet $(\partial E)^{+}\cap\partial
D$. It
suffices to show that $T(d\varphi)=0$. Let $\mu$ denote $H^{n-1}$
restricted to $(\partial E)^{+}$. Appealing to Lemma
\thmref{forms}, we can produce a sequence of Lipschitz functions
$\{\omega_{i}\}$ such that 
$$
\displaylines{
\omega_{i}\to1\;\mu\;\hbox{\rm a.e.}\cr
\abs{\nabla \omega_{i}}\to0\;\mu\;\hbox{\rm a.e.}\cr
\omega_{i}\;\hbox{\rm vanishes in a neighborhood of}\;S\cr
\int_{(\partial E)^{+}}\abs{\nabla \omega_{i}}\;d\mu\to0.\cr}
$$
Thus, we obtain
$$
\eqalign{
0=T(d(\varphi\omega_{i}))
&=T(d\varphi\wedge\omega_{i})+T(\varphi \wedge d\omega_{i})\cr
&=\int_{(\partial E)^{+}}d\varphi\wedge\omega_{i}+\int_{(\partial
E)^{+}}\varphi\wedge d\omega_{i}.\cr}
$$
The first integral tends to
$$
\int_{(\partial E)^{+}}d\varphi=T(d\varphi)
$$
while the second tends to $0$. Thus,
$T(d\varphi)=0$. \qed
\bigskip
Let $E$ denote a 
minimizer of \eqnref{problem}, where $\Omega$ is strictly convex with
$C^2$ boundary. Since $\partial E$
is locally an $n-1$-manifold of class $C^{1}$ except for a singular
set $S$ whose Hausdorff dimension does not exceed $n-8$, it follows
that $\partial E$ can be regarded as an oriented $n-1$ integral
current whose boundary is $0$; i.e. an oriented $n-1$ integral cycle.

Let $T$ denote the $n-1$ integral current represented by $(\partial
E)\cap H^{+}$. Since $\partial E$ is of class $C^{1,1}$ in a
neighborhood of each point of $(\partial E)\cap(\partial \Omega)$, it
follows that the tangent cone to $\partial E$ at such points is in
fact a tangent plane. Consequently, $\partial E$ is analytic near
such points and therefore 
the singular set $S$ of $\partial E$ lies in the
interior of $(\partial E)\cap H^{+}$.
We know from Lemma
\thmref{approxid} that the boundary of $T$
is the $n-2$-sphere determined by $\partial D_{B_{\Omega }}$, the equator 
of $B_{\Omega}$. Let
$p\colon\R^{n}\to T_{B_{\Omega }}$
denote the orthogonal projection and consider the current
$R:=p_{\#}(T)$. Note that $\partial R=p_{\#}(\partial
T)=\partial D_{B_{\Omega}}$. Furthermore, $D_{B_{\Omega}}$ is the
unique current in $T_{B_{\Omega }}$
whose boundary is $\partial D_{B_{\Omega}}$ and therefore, we conclude
that $R=D_{B_{\Omega}}$. Let 
us consider the action of $R$ operating on an $n-1$-form $\varphi$.
For this we will let $\alpha(x)$ denote the Grassman $(n-1)$-vector
of norm one that is in the tangent plane orthogonal to $\nu(E,x)$,
the exterior normal to $E$ at $x$. $\alpha(x)$ is chosen in such a way that
$\alpha(x)\wedge \nu(E,x)$ forms the Grassman unit $n$-vector that
induces a positive orientation of $\R^{n}$. Also, we let
$dp(\alpha(x))$ denote the value of the differential of $p$ operating
on $\alpha(x)$. Then, with the help of \eqnref{intmul}, we have
$$
\eqalign{
R(\varphi)&=T(p^{\#}\varphi)\cr
&=\int_{(\partial E)\cap H^{+}}p^{\#}\varphi\cdot\alpha\cr
&=\int_{(\partial E)\cap H^{+}}\varphi[p(x)]\cdot
dp(\alpha(x))\;dH^{n-1}(x)\cr 
&=\int_{D_{B_{\Omega}}}\varphi(y)[N^{+}(p,\partial E,y)
-N^{-}(p,\partial E,y)]\;dy\cr}
$$
where $N^{+}(p,\partial E,y)$ denotes the number of points of
$p^{-1}(y)\cap\partial E$ at which $Jp$ is positive and similarly, 
$N^{-}(p,\partial E,y)$ denotes the number of points of 
$p^{-1}(y)\cap\partial E$ at which $Jp$ is negative. Since $R=D_{B_{\Omega}}$,
we conclude that $$ N^{+}(p,\partial E,y)-N^{-}(p,\partial
E,y)=1\eqnlbl{sumdeg} $$  for almost all $y\in D_{B_{\Omega}}$.\qed \bigskip
{\bf\thmlbl{HisK} Lemma.} {\it Assume that $\Omega$ is bounded, strictly
convex, has a $C^2$ boundary, and satisfies a great circle condition.
Let $H$ denote the convex hull for any
minimizer $E$ of the variational problem 
{\rm \eqnref{problem}}. Then there is a constant $K$ such that $\mch=K$ at
$H^{n-1}$ almost all points of $(\partial H)\cap\Omega$.}
\medskip
{\bf Proof.} First, we recall that
$\partial E\cap\overline{\Omega}$ is $C^{1}$ at all of its points except
for a singular set $S\subset\partial E\cap\Omega$ whose Hausdorff
dimension does not exceed $n-8$. Furthermore, we know that $\partial
E\cap\Omega$ is real analytic at all points away from $S$ and that
$\partial H$ is $C^{1,1}$. Finally, we know that $E$ contains
$B_{\Omega}$. Let $(\partial E)^{+}$ and $(\partial H)^{+}$ denote the parts
of $\partial E$ and $\partial H$ respectively that lie above the
equatorial plane $P$ 
of $B_{\Omega}$. Let $p\colon\R^{n}\to P$
denote the
orthogonal projection. The mean curvature of $\partial E$ is equal to
a constant $K$ at all points of $\partial E\cap(\Omega-S)$. Let $X$
denote the vertical unit vector. We wish to apply \eqnref{divmc} with
$ (\partial E)^{+}$ replacing $ M$. Referring to the proof of
Lemma \thmref{approxid}, we see that this can be done in spite of
the
singular set $ S\in(\partial E)^{+}$. Thus, applying \eqnref{divmc},
we obtain 
$$
\int_{(\partial H)^{+}}{\cal H}_{\partial H}X\cdot\nu_{H}\;dH^{n-1}
=\int_{(\partial
E)^{+}}{\cal H}_{\partial E}X\cdot\nu_{E}\;dH^{n-1} \eqnlbl{equalint}
$$
where $\nu_{H}$ and $\nu_{E}$ denote the unit exterior normals to $H$
and $E$ respectively.
Let
$$
\displaylines{
A=(\partial E)^{+}\cap(\partial H)^{+}\cr
B=((\partial H)^{+}-A)\cap\{x:{\cal H}_{\partial H}(x)<K\}\cr
C=((\partial H)^{+}-A)\cap\{x:{\cal H}_{\partial H}(x)=K\}.\cr}
$$
Since $\mch \le K$ $H^{n-1}$-a.e. in $(\partial H)^{+}\cap\Omega$, it
suffices to prove 
that 
$$
H^{n-1}(B)=0.\eqnlbl{toprove}
$$
Observe that both $B$ and $C$ are subsets of $\partial H^+$. Note
also that $A,B,$ and $C$ are mutually disjoint subsets of $(\partial H)^{+}$
with $H^{n-1}[(\partial H)^{+}-(A\cup B\cup C)]=0$. Thus, $p(A),p(B)$
and $p(C)$ are mutually disjoint 
and their union occupies almost all of $D_{B_{\Omega}}$.
Clearly, $\nu_{E}$ and $\nu_{H}$ as well as  $\mch$ and
$\mce$ agree $H^{n-1}$ almost everywhere on $A$. Therefore,
$$
\int_{A}\mch X\cdot\nu_{H}\;dH^{n-1}=\int_{A}\mce
X\cdot\nu_{E}\;dH^{n-1}\;.\eqnlbl{ona} 
$$
Since $X\cdot\nu_{H}$ is the Jacobian of the mapping $p\colon\partial
H^{+}\to D_{B_{\Omega}}$, it follows from \eqnref{intmul} that
$$\displaylines{
\int_{B}\mch X\cdot\nu_{H}\;dH^{n-1}<K H^{n-1}[p(B)],\cr
\int_{C}\mch X\cdot\nu_{H}\;dH^{n-1}=K H^{n-1}[p(C)].\cr}
$$
Now let 
$$
\eqalign{
A^{*}&=((\partial E)^{+})\cap p^{-1}[p(A)],\cr
B^{*}&=((\partial E)^{+})\cap p^{-1}[p(B)],\cr
C^{*}&=((\partial E)^{+})\cap p^{-1}[p(C)].\cr}
$$

Next, observe that both $B^{*}$ and $C^{*}$ are subsets of $\Omega$.
To see this, consider $x\in B^{*}$. If it were true that
$x\in B^{*}\cap\partial \Omega $, then $x\in(\partial H)^{+}$ and
thus $x\in A$. This is impossible since $p(A)$ and $p(B)$ are
disjoint. A similar argument holds
for $C^{*}$.
Referring to \eqnref{intmul} and \eqnref{sumdeg}, we obtain
$$
\eqalign{
&\int_{B^{*}}\mce X\cdot\nu_{E}\;dH^{n-1}\cr
&\quad=K\int_{B^{*}\cap\{x:X\cdot\nu_{E}(x)>0\}} X\cdot\nu_{E}\;dH^{n-1}
+K\int_{B^{*}\cap\{x:X\cdot\nu_{E}(x)<0\}}
X\cdot\nu_{E}\;dH^{n-1}\cr 
&\quad=K\int_{p(B^{*})}N^{+}(p,\partial E,y)-N^{-}(p,\partial
E,y)\;dH^{n-1}(y)\cr 
&\quad=K H^{n-1}[p(B^{*})]\cr
&\quad=K H^{n-1}[p(B)].\cr}
$$
Similarly, 
$$
\int_{C^{*}}\mce X\cdot\nu_{E}\;dH^{n-1}
=K H^{n-1}[p(C^{*})]=K H^{n-1}[p(C)]
$$
and
$$ 
\int_{A^*} K X\cdot\nu_E\;dH^{n-1}=K H^{n-1}(p(A)).
$$
Finally, in view of the fact that $A\subset(\partial H)^{+}$ and
therefore that $N^{+}(p,A,y)=1$ and $N^{-}(p,A,y)=0$ for
$H^{n-1}$-almost all $y\in p(A)$, we obtain
$$
\int_{A} K X\cdot\nu_E\;dH^{n-1}=
K H^{n-1}(p(A)).
$$
Now, using the facts that $A^{*}-A\subset\Omega$ and 
$\mce=K$ on $A^{*}-A-S$, we obtain
$$
\eqalign{
\int_{A^*}&\mce X\cdot\nu_{E}\;dH^{n-1}\cr
&\qquad=\int_{A^*}K X\cdot\nu_{E}\;dH^{n-1}
+\int_{A^*}(\mce-K) X\cdot\nu_{E}\;dH^{n-1}\cr
&\qquad=\int_{A^*}K X\cdot\nu_{E}\;dH^{n-1} +\int_{A}(\mce-K)
X\cdot\nu_{E}\;dH^{n-1} \cr 
&\qquad=K H^{n-1}(p(A))-K H^{n-1}(p(A))+\int_{A}\mce
X\cdot\nu_{E}\;dH^{n-1}\cr 
&\qquad=\int_{A}\mce X\cdot\nu_{E}\;dH^{n-1}.\cr}
$$
Under the assumption $H^{n-1}(B)>0$, we would obtain
$$
\eqalign{
\int_{(\partial H)^{+}}\mch X\cdot\nu_{H}\;dH^{n-1}&
<\int_{A}\mch X\cdot\nu_{H}\;dH^{n-1} 
+K H^{n-1}[p(B)]+K H^{n-1}[p(C)]\cr
&=\int_{A}\mce X\cdot\nu_{E}\;dH^{n-1}+K H^{n-1}[p(B^{*})]+K
H^{n-1}[p(C^{*})]\cr 
&=\int_{A^{*}}\mce X\cdot\nu_{E}\;dH^{n-1}+K H^{n-1}[p(B^{*})]+K
H^{n-1}[p(C^{*})]\cr 
&=\int_{A^{*}}\mce X\cdot\nu_{E}\;dH^{n-1}+\int_{B^{*}}\mce
X\cdot\nu_{E}\;dH^{n-1}\cr
&\qquad+\int_{C^{*}}\mce X\cdot\nu_{E}\;dH^{n-1}\cr 
&=\int_{A^{*}\cup B^{*}\cup C^{*}}\mce X\cdot\nu_{E}\;dH^{n-1}\cr
&\leq\int_{(\partial E)^{+}}\mce X\cdot\nu_{E}\;dH^{n-1},\cr}
$$
where we have used that $A^{*},B^{*}$ and $C^{*}$ are mutually disjoint.
This would contradict \eqnref{equalint}, thus establishing
\eqnref{toprove}. \qed 
\bigskip
A function $u\in C^{1}(W)$ is called a {\it weak subsolution
{\rm(}supersolution{\rm )} of the equation of constant $K$ mean
curvature}
if
$$
Mu(\varphi)=\int_{W}\frac{\nabla u\cdot\nabla
\varphi}{\sqrt{1+\abs{\nabla u}^{2}}}-K\varphi\,dx\leq0\quad(\geq 0)
$$ 
whenever $\varphi\in C^{1}_{0}(W),\; \varphi \geq 0$.

We note that if $u\in C^{1,1}$ and classically satisfies the equation
of constant mean curvature equation almost everywhere, then $u$ is a
weak solution. 

The following result will be stated in the context of $R^{n-1}$
because of its applications in the subsequent
development. 
\bigskip
{\bf \thmlbl{har} Lemma.} {Suppose $W$ is an open subset of
$R^{n-1}$. If $u_{1}, u_{2}\in C^{1}(W)$ are respectively weak super
and subsolutions of the equation of constant mean curvature in $W$ and if
$u_{1}(x_{0})=u_{2}(x_{0})$ for some $x_{0}\in W$ while $u_{1}(x)\geq
u_{2}(x)$ for all $x\in W$, then
$$
u_{1}(x)=u_{2}(x)
$$
for all $x$ in some closed ball contained in $W$ centered at
$x_{0}$.} 
\medskip
{\bf Proof.} Define
$$
\eqalign{
u_{t}&=tu_{1}+(1-t)u_{2}\;\hbox{\rm for}\;t\in[0,1],\cr
w&=u_{1}-u_{2},\cr
a^{ij}(x)&=\int_{0}^{1}D_{u_{x_{j}}}\left(\frac{D_{i}u_{t}(x)}
{\sqrt{1+\abs{\nabla
u_{t}}^{2}}}\right)\,dt\cr 
&=\int_{0}^{1}\frac{1}{\sqrt{1+\abs{\nabla
u_{t}}^{2}}}\left(\delta_{ij}- 
\frac{D_{i}u_{t}(x)D_{j}u_{t}(x)}{(1+\abs{\nabla
u_{t}}^{2})}\right)\;dt.\cr}
$$
Since both $u_{1}$ and $u_{2}$ are continuously differentiable in
$W$, for each open set $V\subset\subset W$ containing $x_{0}$ there
exists $M>0$ such that $\abs{\nabla u_{t}(x)}\leq M$ for all $x\in V$
and all $t\in[0,1]$. Hence,
$$
\eqalign{
a^{ij}(x)\xi_{i}\xi _{j}&\geq\frac{1}{(1+M^{2})^{1/2}}\abs{\xi
}^{2},\;\hbox{\rm for all}\;\xi \in R^{n-1},x\in V,\cr
\sum_{i,j}a^{ij}(x)^{2}&\leq C,\;\hbox{\rm for all}\;x\in V.\cr
}
$$
For $\varphi\in
C^{1}_{0}(W),\>\varphi\geq0$, we have
$$
\eqalign{
0&\leq Mu_{1}(\varphi)-Mu_{2}(\varphi)\cr
&=\int_{W}\!\!\int_{0}^{1}\frac{d}{dt}\left(\frac{\nabla
u_{t}(x)\cdot\nabla \varphi(x)}{\sqrt{1+\abs{\nabla
u_{t}}^{2}}}\right)\,dt\,dx\cr 
&=\int_{W}a^{ij}(x)D_{j}w(x)D_{i}\varphi(x)\,dx.\cr}
$$
Thus, $w$ is a weak supersolution of the equation
$$
D_{i}(a^{ij}D_{j}w)=0
$$
and since $w\geq0$, the weak Harnack inequality [GT, Theorem 8.18]
yields $$
\left( r^{-n}\int_{B(x_{0},2r)}\abs{w(x)}^{p}\,dx\right)^{1/p}\leq
C\inf_{B(x_{0},r)}w=0 
$$
whenever $1\leq n<n/(n-2)$ and $B(x_{0},4r)\subset W$.\qed
\bigskip

{\bf \thmlbl{econvex} Theorem.} {\it Suppose $\Omega$ is a bounded,
strictly 
convex domain with $C^{2}$ boundary that satisfies a great circle
condition. Then any minimizer $E$ of
the variational problem {\rm \eqnref{problem}} is convex.}
\medskip
{\bf Remark.} Later we show that
neither smoothness of  $\partial\Omega$ nor strict convexity
are required. In addition,
the great circle condition is unnecessary in $\R^{2}$. The same
applies to the uniqueness result below.
 
{\bf Proof.} 
It suffices to show that $ H=E$ where $ H$ denotes the
convex
hull of $ E$. Assume $\partial H\not\subset\partial E$ so there exists
$x\in \partial H\setminus\partial E$.  Thus, as in the proof
of the mean curvature inequality in Theorem \thmref{bhregular}, we
see that $x$ lies in the convex hull $M$ of distinct points
$p_i\in \partial H\cap\partial E$, $i=1,\dots,k,\, k>1$.
Futhermore each $p_{i}$ is an element of $\Omega $ due to the fact
that they all lie in a single support plane
of $H$; hence if one $p_i$ where to lie in $\partial \Omega$ then
they all would, thus contradicting strict convexity.
Referring to 
Lemma \thmref{har}, we see that $ \partial H$ and $ \partial E$ agree
in a neighborhood of the points $p_i$. Since $ M$ is connected,
it follows again from Lemma \thmref{har} that 
$ M\subset\partial E\cap\partial H$, which
contradicts $ x\not\in\partial E$. 
Consequently $\partial H\subset\partial E$ and thus $P(H)\le P(E)$.
However $E\subset H$ so $|E|\le|H|$.  Assume $|E|<|H|$.
Dilate $H$ to obtain $\tilde H\subset \Omega$ satisfying $|\tilde
H|=|E|$.
But then $P(\tilde H)<P(H)\le P(E)$ which contradicts the minimality
of $E$. Thus $|E|=|H|$ so that $E$
and $H$ have the same measure theoretic closure. Hence, due to our
convention concerning 
distinguished representatives for sets of finite perimeter, 
$E=H$ and $E$ is convex.\qed
\bigskip
{\bf\thmlbl{enested} Theorem.} {\it If $\Omega$ 
is as in Theorem {\rm \thmref{econvex}} 
then perimeter minimizers with measure exceeding $|B_\Omega|$ are
nested and unique. That is, if $E$ and $F$ are perimeter minimizers
then
$$
|B_\Omega|\le|F|<|E| \Longrightarrow E\subset F   \eqnlbl{nested}
$$
and 
$$
|B_\Omega|\le |F| = |E| \Longrightarrow F=E      \eqnlbl{eunique}
$$
In addition perimeter minimizers have disjoint boundaries 
relative to $\Omega$ in the sense
that 
$$
|B_\Omega|\le|F|<|E|  \Longrightarrow
\partial F\cap\partial E\subset\partial\Omega.
$$}
\smallskip
{\bf Remark.} Note that the assumption of convexity can be relaxed.
It is only required that the intersection of $\Omega$ with any
vertical line is an 
interval. In addition $\partial\Omega$ must not contain vertical line
segments.
\medskip
{\bf Proof.}  To prove \eqnref{nested} we 
argue by contradiction. If $E$ and $F$ are perimeter minimizers satisfying
$|B_\Omega|\le|F|<|E|$ assume $F$ is not a subset of $E$.
From Theorems \thmref{ballinmin} and \thmref{econvex} we see that $E$
and
$F$
are convex and contain $B_\Omega$. 
Since $F$ is not a subset of $E$ one can employ the proof of
\thmref{ballinmin}, 
with $F$ playing the role of $B_\Omega$, to prove that 
there is a second perimeter minimizer $E^{*}$ 
which contains $F$ and satisfies $|E^{*}| = |E|$. Let $H$ be the analog
of $H$ in the proof of Theorem \thmref{ballinmin} and let $D^{\circ}$
denote the interior of $D:=D_{B_{\Omega }}$.

We will use the  properties of perimeter minimizers to show that
$\partial H$ and $\partial(H\cup E)$ are analytic and coincide on
some open set. By connectedness, this will show they are
identical, thus establishing the desired contradiction.

Let $O$ be the interior of $\partial H\setminus (H\cup E)^\circ$ 
relative to $\partial H$, and $\partial O$ represent the boundary of $O$
relative to $\partial H$.
Assume there exists a point
$$
x\in\partial O\cap p^{-1}(D^\circ).
$$
Note that $x\in \partial H\cap \partial E\cap \partial (H\cup E)\cap
p^{-1}(D^\circ)$. Let $y$ be the point on $\partial F\cap
\partial E^{*}$ which was translated (as in the definition of $H$) to
$x$. Since $\partial O$ has positive $H^{n-2}$ measure ($\partial
O\cap p^{-1}(D^\circ)\ne\nullset$) we can assume $y\notin S,\,S$ being
the singular set for $E^{*}$.

Since $x\in\partial E\subset\overline\Omega,\, y$ lies in $\Omega$ and
consequently $\partial H$ is analytic in a neighborhood of $x$ since
$\partial F$ is analytic in a neighborhood of $y$. Similarly  $H\cup
E$ inherits analyticity (in a neighborhood of $x$) from
$\partial E^{*}$ since $x\in \partial(H\cup E)$ and $y\notin S$.
However $\partial H\cap O\subset \partial (H\cup E)$ so $\partial H$
and $\partial (H\cup E)$ coincide on open (relative to $\partial H$)
subsets of any neighborhood of $x$ so by analyticity $\partial H$
coincides with $\partial (H\cup E)$ in some  neighborhood of $x$. But
this contradicts $x\in\partial O$ so $\partial O\cap p^{-1}(D^\circ)$
is empty.

Note that $(\partial H\setminus E)\cap p^{-1}(D^\circ)$ contains
points lying both  above and below $D$ since $\Omega$ is strictly
convex and $H$ is the hull of the translated halves of $F$ (which
contain the hemispheres of the largest ball $B_\Omega$). Thus the
same is true of  $O\cap p^{-1}(D^\circ)$. Combined with $\partial
O\cap p^{-1}(D^\circ)=\nullset$, this implies $\partial H\cap
p^{-1}(D^\circ)\cap E^\circ=\nullset$, i.e. $E\subset H$. Of course
this is absurd since $|E|>|F|=|H|$. Thus  the assumption that $F$ is
not contained in $E$ is false i.e. $F\subset E$ as required.

Now assume  that $|B_\Omega|<|F| = |E|=v$. Choosing a sequence of
perimeter minimizers $F_i$ of measure $v_i\uparrow v$, it follows
from \eqnref{nested} that  $F_i\subset E\cap F$.  Consequently
$|E\cap F|=v$ and so $E=F$.

To prove that minimizers are strictly nested in the sense defined
above assume that $|B_\Omega|\le|F|<|E|$ and so $F\subset E$. Assume
in addition that $G:=(\partial F\cap\partial E\cap\Omega)^\circ$ is
not empty. Since $F,\, E$ are analytic  in $\Omega$ and nested, it is
clear that $\mcf\ge\mce$ at points in $G$. Given that $\mcf,\,\mce$
are  constants, say $k_f,\, k_e$, in $\Omega$ and equal almost
everywhere on $\partial F\cap\partial E\cap\partial\Omega$, we may
derive a contradiction from $k_f\ge k_e$ through the use of
\eqnref{divmc}. In fact, we obtain 
$$
\int_D \mcf'\;dH^{n-1}=H^{n-2}(\partial D)=
\int_D \mce'\;dH^{n-1} \eqnlbl{sameint}
$$
where $\mcf'(x):=\mcf(p^{-1}(x)\cap\partial F)$ and
$\mce'(x):=\mce(p^{-1}(x)\cap\partial E)$.
However with
$A:=\partial F\cap\partial \Omega$ and $B:=\partial F\cap\Omega$,
we see that
$$
\eqalign{
\int_D \mce'\;dH^{n-1}
&=\int_{p(A)} \mcf'\;dH^{n-1}+\int_{p(B)}\mce'\;dH^{n-1}\cr
&\le \int_{p(A)}\mcf'\;dH^{n-1}+\int_{p(B)}\mcf'\;dH^{n-1}\cr
&=\int_D\mcf'\;dH^{n-1},\cr}
\eqnlbl{equality}
$$
and thus we have equality due to \eqnref{sameint}. Therefore, 
$\mce'=k_e=k_f=\mcf'$
on $p(B)$.
However, since $\mce'=\mcf'$ almost everywhere on $p(A)$, we obtain
$$
\mce'=k_e=k_f=\mcf' \qquad \hbox{$H^{n-1}$-almost everywhere on} D.
$$
Thus for $x\in \partial F\cap\partial E$, apply Lemma \thmref{har} to
conclude 
that $\partial F$ and $\partial E$ coincide in a neighborhood of $x$.
Thus $p(\partial F\cap\partial E)$ is both open and closed relative
to $D^\circ$, and therefore contains $D^\circ$, a
contradiction
since $\abs{F}<\abs{E}$.\qed

\bigskip

We now dispense with the assumptions of strict convexity and
smoothness of  $\partial\Omega$.  When the assumption of strict
convexity is dropped, complications arise because there is no longer
a unique largest ball in $\Omega$. Eliminating  the smoothness
assumption on the boundary forces us to take limits of perimeter
minimizers, and to establish convexity of all perimeter minimizers
through a uniqueness theorem.

One interesting observation is that a perimeter minimizer can be
thought of as a smooth approximation of $\Omega$, especially when its
measure is close to that of $\Omega$. This is due to the fact that
even after we have dispensed with the smoothness assumption on
$\partial \Omega$ perimeter minimizers still have $C^{1,1}$
boundaries.

For the proof of Theorem \thmref{notstrict} below, we need the
following lemma.
\bigskip
{\bf \thmlbl{smoothconvex} Lemma.} {\it Let $ a<c<b$ and let 
$I_{1},I_{2}$ denote the closed intervals $ [a,c]$ and $ [c,b]$,
respectively. Let $ f_{1}$ and $ f_{2}$ be functions such that 
$ f_{i}\in C^{2}(I_{i}),\;i=1,2$, 
with $ f_{1}(c)=f_{2}(c)$. Furthermore, assume there are constants
$ c_{1},c_{2}$ and $ c_{3}$ such that 
\medskip
\itemitem{\rm(i)} $f_{i}''\leq c_{1}<0$ on $ I_{i},\;i=1,2$,
\medskip
\itemitem{\rm(ii)} $f_{1}'\geq c_{2}>0$ on $ I_{1}$ and 
$ f_{2}'\leq c_{3}<0$ on $ I_{2}$.
\medskip\noindent
Then, there exists a $C^{2}$, strictly concave function $ g$ on 
$[a,b]$ such that $ g$ is uniformly close to $ f$ on $[a,b]$ and that
$ g=f$ on the complement of any given open interval containing $ c$.}
\medskip
{\bf Proof.} A given open interval containing $c$ in turn contains
 an open interval $ I=(a',b')$ with $ c\in I$
determined by the constants $ c_{1},c_{2}$ and $ c_{3}$ such that the
following three conditions hold:
\medskip
\itemitem{(i)} There are points $ x_{1},x_{2}\in I$ with
$x_{1}<c<x_{2}$ such that 
$ f_{1}(x_{1})=f_{2}(x_{2})$. 
\medskip
\itemitem{(ii)} There are polynomials $ p_{i}$ of
degree $ 2$ ($i$=1,2) such that $p_{i}(x_{i})=f_{i}(x_{i}) $ and such
that the functions 
$$
h_{1}(x):=\cases{
			f_{1}(x)& for $ a\leq x\leq x_{1}$\cr
			p_{1}(x)& for $ x_{1}\leq x\leq c$\cr}
\qquad h_{2}(x):= \cases{
			f_{2}(x)& for $ x_{2}\leq x\leq b$\cr
			p_{2}(x)& for $ c\leq x\leq x_{2}$\cr}
$$
are $ C^{2}$ and strictly concave on $ I_{i}$.
\medskip
\itemitem{(iii)} There is a point $c'\in I $ such that
$ h_{1}(c')=h_{2}(c')$. 
\medskip\noindent
Thus, the function 
$$
h:=\cases{
h_{1}& on $ [a,c']$\cr
h_{2}& on $ [c',b]$\cr}
$$
is strictly concave on $ [a,b]$. We now will mollify $ h$ restricted
to $ I$ by using a smooth 
mollifying kernel $ \varphi $ with the property that 
$$
\varphi _{\varepsilon }\ast p(x)=p(x) 
$$
whenever $ p$ is a polynomial of degree 2, $\varepsilon >0$, and
$ x\in\R$, cf. [Z, Lemma 3.5.6]. Thus, for sufficiently small
$ \varepsilon >0,\;\varphi _{\varepsilon }\ast h(x)=h(x)$, for
$ x\in (a'+\varepsilon ,c'-\varepsilon )
\cup(c'+\varepsilon,b'-\varepsilon )$. Also, 
$ \varphi _{\varepsilon }\ast h$ is strictly concave since $ h$ is.
Thus, our desired function $ g$ is defined by 
$$
g(x) =\cases{
	    h(x)& for $ a\leq x\leq a'+\varepsilon $\cr
	    \varphi _{\varepsilon }\ast h(x)& for 
	    $ a'+\varepsilon < x <c'-\varepsilon $\cr
	    h(x)& for $ c'-\varepsilon \leq x\leq b$.\cr}
$$
\qed

We define $H_\Omega$ to be the union of all largest balls in
$\Omega$.  Thus $H_\Omega$ is the convex hull of the two largest
balls which are furthest apart. $H_\Omega$ essentially plays the role
of $B_{\Omega }$. 
\bigskip 
{\bf \thmlbl{notstrict} Theorem.} {\it
Suppose $\Omega$ is a bounded, convex domain  that satisfies a great
circle condition. Given $v,\,|H_\Omega|\le v<|\Omega|$ there is a
unique  minimizer $E$ with $|E|=v$ of the variational problem {\rm
\eqnref{problem}}. $E$ is convex with $C^{1,1}$ boundary. Such
minimizers are nested with disjoint boundaries relative to $\Omega$
 as in Theorem \thmref{enested}. 
If $|B_\Omega|< v \le|H_\Omega|$  then any minimizer $E$ is the
convex hull of two largest balls (clearly uniqueness is lost
for $v<|H_\Omega|$).}
\medskip 
{\bf Proof.} We first smooth $\Omega$ and then establish the
existence of a nested family of convex perimeter minimizers by taking
limits.  We finish by adapting the uniqueness result of Theorem
\thmref{enested} and the proof of disjointness of boundaries.

Let $T_{B_\Omega}$ be the hyperplane which
intersects orthogonally the midpoint
of the line segment joining the centers of the two largest balls
whose hull forms $H_\Omega$. Think of the  ``vertical'' axis as
coinciding with this line segment and take the origin of our
coordinate system to be the midpoint just
 mentioned.  As defined
previously $p$ is orthogonal  projection 
onto $T_{B_\Omega}$. Let $B_\Omega$ be
the largest ball in $\Omega$ with equatorial plane in $T_{B_\Omega}$.   Let
$D_{B_\Omega}=p(\Omega)$ so $D_{B_\Omega}$ 
is an $(n-1)$-ball.  Let $C$ be the interior
of the union of a closed 
 right circular cone with base $D_{B_\Omega}$
with its reflection across $T_{B_\Omega}$.  Let $B$ be the largest ball in $C$
and note that $C\setminus B$ has three components (four in $\R^2$).
Let $C_0$ denote the component (or union of two components in $\R^2$)
which intersects $D_{B_\Omega}$ and consider the
set $C_1=C\setminus \overline{C_0}$.

First we show that $\Omega$ can be approximated arbitrarily closely
by strictly convex sets satisfying a great circle condition, then we 
will approximate the later by sets with $ C^{2}$ boundary of the same
type. Note that $\Omega\cap C_1$ is convex and satisfies a great
circle condition with $B$ being the largest ball. Also
$\partial(\Omega\cap C_1)$ consists of the union of the graphs of
functions $f_i,\,i=1,2,\, f_1\ge0,\,f_2\le 0$. Let $\Omega'$ be the
set whose boundary is the union of the graphs of  $f_1+\varepsilon b,\,
f_2-\varepsilon b$ where $\varepsilon>0$ and  $b$ is the function
whose graph is the upper hemishere of $B$. Note that $\Omega'$ is
strictly convex and satisfies a great circle condition. Also, as
$\varepsilon\to 0$, $C$ approaches a cylinder, and $\Omega'\to\Omega$
in the Hausdorff sense.

We now may 
assume with out loss of generality that $\Omega$ is strictly convex.
Consider $G=\Omega \cap C$. Note that $\partial G$ is the union of
graphs of 
$f_i:\overline D_{B_\Omega}\rightarrow \R,\, i=1,2$ with $f_1\ge0,\, f_2\le
0$. Given $r>0$ let $B_r$ be the  ball of radius $r$ concentric to $B_\Omega$,
$D_r=D_{B_\Omega}\cap B_r$,and $R$ the radius of $B_\Omega$. Also let
$\bar r$ be the distance from $\partial
B_\Omega\cap\partial C$ to the vertical axis.

Consider $\varepsilon, \, 0<\varepsilon<<R$. For a smooth
radially symmetric approximate identity $\eta_\varepsilon$ supported
in $B_\varepsilon$  let $f_\varepsilon=f_1*\eta_\varepsilon$. Thus
$f_\varepsilon$ is defined in $D_{R-\varepsilon}$ and is a surface of
revolution in  $A_\varepsilon=D_{R-\varepsilon}\setminus D_{\bar
r+\varepsilon}$ 

Now consider $\delta>0$ such that $\bar r<R-\delta$ but 
$\partial B_{R-\delta}$ does not intersect $\partial C$.
Take $\varepsilon$ small enough that the graph of $f_\varepsilon$
does not intersect $\partial B_{R-\delta}$.  Let $g_\varepsilon:[\bar
r+\varepsilon, R-\varepsilon]\rightarrow \R$ be the function the
rotation of whose graph around the vertical axis produces the graph of
$f_\varepsilon$ over $A_\varepsilon$. In the $r,\,z$ plane let $C_2$
be a circle of radius $s>>R$ with center on the  negative $r$ axis
which passes through $(R-\delta,0)$. Let $c:[\bar r+\varepsilon,
R-\varepsilon]\rightarrow \R$ be the function whose graph lies in
the upper half of $C_2$ and define
$h_\varepsilon=\min(g_\varepsilon,\,c)$ on $[\bar r+\varepsilon,\,R-\delta]$.
Note $h_\varepsilon$ is a
strictly concave function and is smooth except at the point $q$ of
intersection of the graphs of $g_\varepsilon$ and $c$ (which exists if 
$s$ is large enough). Now employ
Lemma \thmref{smoothconvex} to alter 
$h_{\varepsilon }$ in a small neighborhood of $q$ to produce a $C^2$
function which is still strictly concave.

Consider the surface obtained by taking the union of the surface of
revolution formed by rotating the graph of the smoothed
$h_\varepsilon$ with the graph of $f_\varepsilon$ over $D_{\bar
r+\varepsilon}$. This is a $C^2$ surface and when combined with a
similarly constructed surface for $f_2$ produces the boundary of a
strictly convex set $\Omega_\varepsilon$. Note that $\partial
\Omega_\varepsilon$ is $C^2$ and that $\Omega_\varepsilon$ satisfies a
great circle condition with $B_{R-\delta}$ being the largest ball.
Also as $C$ approaches a cylinder and $\delta, \varepsilon\to 0$ we have
$\Omega_\varepsilon\to\Omega$ in the Hausdorff sense as required. To
make the process of taking limits easier in
the following we can dilate the sets
$\Omega_\varepsilon$ a small amount so they contain $\Omega$.

Thus there exists a sequence of $C^2$ strictly convex sets $\Omega_n$
which contain $\Omega$, satisfy a great circle condition, and which
converge to $\Omega$ in the Hausdorff sense. For $v,\,|B_\Omega|<
v\le|\Omega|$ (and $n$ large enough so $|B_{\Omega_n}|<v$) let
$E_n(v)$   be the unique perimeter minimizer in $\Omega_n$ of measure
$v$. It is easy to see that for a dense set of $v_i$'s we can, by
repeatedly extracting subsequences and diagonalizing, construct a
subsequence of $E_n$ such that for all $i$, $E_n(v_i)$ converges (on
the subsequence) to   $E(v_i)$, a subset of $\Omega$, in the Hausdorff
sense. Nestedness and convexity are clearly inherited. Thus  taking
intersections of appropriate $E(v_i)$ we extend the definition of
$E(v)$ to  all $v,\,|B_\Omega|< v<|\Omega|$.  Nestedness allows us
to extend convergence to all such $v$.

We claim that the sets $E(v)$ are perimeter minimizers relative to
$\Omega$.  To see this note that given any set $F\subset\Omega$ with
$|F|=v$ we have $F\subset\Omega_n$ since $\Omega\subset\Omega_n$;
consequently by lower semicontinuity  of perimeter we have
$$
P(E(v))\le \liminf P(E_{n}(v)) \le P(F)
$$
(with the liminf taken over the subsequence) i.e. $E(v)$ is a
perimeter minimizer. 

For $v,\,0\le v\le|H_\Omega|$  we can characterize perimeter
minimizers. Assume $E$ is a perimeter minimizer of measure $v$. If
$0<v\le|B_\Omega|$ then $E$ is clearly a ball. If $|B_\Omega|<
v\le|H_\Omega|$ we claim that $E$ is the convex hull of two largest
balls in $\overline\Omega$. In proving this we will 
also prove for $v \ge |H_\Omega|$  that	 any perimeter minimizer $E$
satisfies $H_\Omega\subset E$. Assume $|B_\Omega|< v$.
Consider the following extension of the proof of Theorem
\thmref{ballinmin}. As it stands the proof of Theorem
\thmref{ballinmin} implies that $E$  contains a largest (in $\overline\Omega$)
ball. In fact one can conclude much more. Let $B_1,\,B_2$ be the
closed balls whose convex hull is $H_\Omega$, let $\ell$ be the line
through their centers, and  consider any set $H$ which is the convex
hull of two translates of $B_1$ with centers on $\ell$ such that
$|H\cap E|=|B_1|$ and $H\cap H_\Omega$ contains a translate  of $B_1$.
A mild variation in the proof of Theorem \thmref{ballinmin} shows that
$H\cap E$ is a translate of $B_1$. We claim that this implies that $E\cap
H_\Omega$ is the convex hull of two translates of $B_1$.  To see this let
$B_3,\,B_4\subset E$ be distinct translates 
of $B_1$ with $x$ being the midpoint 
between their centers.  Since the hull of $B_3,\,B_4$ has measure
larger than $|B_1|$ construct $H$ as above using translates of $B_1$ 
placed symmetrically with respect to $x$. However $H\cap E$ is a 
translate of $B_1$.  Thus there is a translate of $B_1$ 
contained in $E$ lying strictly between ant two such balls.
Therfore the centers of such balls form an interval in $\ell$

Now take $\ell$ to be the vertical axis with $B_1$ lying
above $B_2$, let $B_u,\,B_l$ be the the uppermost and lowest 
translates of $B_1$ in $E$, and $E_u,\,E_l$ the parts of $E$
strictly above and below $B_u,\,B_l$ respectively.
Assume $E_u$ is not empty so $|E_u|\ne 0$. If $B_u\ne B_1$ 
construct $H$ as above by translating hemispheres of $B_1$
so that $H$ contains subsets of positive measure from both $E\cap
H_\Omega$ and $E_u$.  However this is a contradiction since by the
above $E\cap H$ is  a translate of $B_1$ which cannot possibly
intersect $E_u$. Thus $E_u$ not empty implies $B_u=B_1$.
Similarly $E_l$ not empty implies $B_l=B_2$.
This establishes the claim.

Moreover one can conclude that $v\ge|H_\Omega|$ implies $H_\Omega\subset E$.
To see this note if $v\ge|H_\Omega|$ then at least one of $E_u,\,E_l$ is 
nonempty. If both are nonempty then $H_\Omega\subset E$ as claimed.
If only one is nonempty, say $E_u$, then translate $E$ as
far down as possible while remaining in $\overline\Omega$ to form a set $E^*$
which contains $B_2$ ($E_l$ is empty). 
Note that $E^*$ is also a perimeter minimizer of
measure $v$. Thus $E_u^*$ nonempty, i.e. 
$E^*=E$ with $H_\Omega\subset E$ as required.

Now that we have characterized perimeter minimizers for $v,\,0<v\le |H_\Omega|$
we can redefine $E(v)$ so that $E(v)$ is the convex hull of 
two translates of $B_1$, symmetrically 
placed in $H_\Omega$, if $|B_1|<v\le |H_\Omega|$, 
and  $E(v)$ is a symmetrically placed ball if $0<v\le |B_1|$.
Thus  we have a nested collection of 
convex perimeter minimizers which can be used to establish
uniqueness. Given $\bar v$, 
$|H_\Omega|<\bar v$, assume that $E$ is a perimeter minimizer with
measure $\bar v$. Recall from  above  that $H_\Omega\subset E$.
Before proceeding we define an auxiliary collection 
$\{H(v):\, |H_\Omega|\le v\le \bar v \}$, $H(v)$ defined analogously to $H$ in
Theorem \thmref{ballinmin} by translating the halves of $E(v)$ 
the least possible amount such that the resultant hull $H(v)$
satisfies $|H(v)\cap E|=v$. Note that the sets $H(v)$ are
nested since if $|H_\Omega|\le v<w$ and one translates the halves of
$E(w)$ the same distance as for $E(v)$ in the definition of $H(v)$, and 
calls the hull of the translated halves $\tilde H$ then
$|\tilde H\setminus H(v)|=w-v$ so
$|\tilde H\cap E|=|(\tilde H\setminus H(v))\cap E|
+|H(v)\cap E|\le (w-v)+v=w$ i.e. $H(v) \subset \tilde H\subset H(w)$
as required.

Let $v_0=\sup\{v:E(v)\subset E\}$. If $v_0=|E|$ then $E=E(v_0)$,
otherwise $v_0<|E|$ so $E^\circ\setminus E(v_0)$
is not empty. Let $B$ be a closed ball 
of positive radius in $E^\circ\setminus E(v_0)$,
$v_1=\sup\{v:H(v)\cap B \hbox{ is empty }\}$, and 
$v_2=\inf\{v:B\subset H(v)\}$. Clearly 
$v_2=|H(v_2)\cap E|\ge |H(v_1)\cap E|+|B|=v_1+|B|$
so choosing $v,\,v_1<v<v_2$ we see that $B$ contains points in $H(v)$
and its complement. Consequently $\partial H(v)$ intersects $B$.
One can now proceed as in the proof of Theorem \thmref{enested}
with $H$ replaced by $H(v)$ with the following modifications.
In proving that 
$(\partial H(v)\setminus E)\cap p^{-1}(D_{B_\Omega}^\circ)$ is not 
empty one uses the fact proved above that $H_\Omega\subset F$ so that $H(v)$
contains a convex hull of ``largest balls'' which is larger than $H_\Omega$
and thus must intersect 
the complement of $\Omega$. Finally we see that the conclusion
$\partial H(v)\cap p^{-1}(D_{B_\Omega}^\circ)\cap E^\circ=\nullset$ 
is absurd due to our construction
in which $\partial H(v)$ intersects $E^\circ$.  Thus the assumption
that $v_0<|E|$ must be false and consequently $E=E(|E|)$ as
required. 

It remains only to prove the
disjointness result. The proof is identical to that
in Theorem \thmref{enested} 
once we have established that minimizers have $C^{1,1}$
boundaries and satisfy the same mean curvature properties as before.
Assume $|H_\Omega|\le v<|\Omega|$. Let $E_n(v)$ be as above and note that
$k_n$, the constant mean curvature 
associated with $\partial E_n(v)$, is bounded
uniformly in $n$ since as in the proof of Theorem \thmref{enested} we have
$$
H^{n-1}(\partial D_{B_\Omega})=\int_{D_{B_\Omega}} \mcen' 
\ge k_n H^{n-1}(p(\partial E_n(v)\cap\Omega^\circ))
$$
where $H^{n-1}(p(\partial E_n(v)\cap\Omega^\circ))$ 
is uniformly bounded from zero
(on a subsequence) for
geometrical reasons since $E_n(v)$ is convex, 
contains $B_\Omega$, and converges
(on a subsequence) to $E(v)$.  Consequently 
$0\le \mcen \le k_n\le M$ almost everywhere
and we see that $\partial E_n(v)$ is uniformly 
$C^{1,1}$ from which we see that $E(v)$
is $C^{1,1}$ as well. Note that tangent 
planes converge almost everywhere so that
locally first derivatives converge almost everywhere and consequently one can
take limits in the weak definition of mean 
curvature to show that if $k_n\rightarrow k$
(on a subsequence) then $\partial E(v)$ has mean curvature $k$ in the interior
of $\Omega$, and that $\mce\le k$ as required.\qed

\bigskip

{\bf\thmlbl{2d} Theorem.} {\it If $n=2$, and $\Omega$ is as in Theorem
\thmref{notstrict}, except that the great circle condition is not
assumed, then the results of  Theorem \thmref{notstrict} still hold.
Furthermore, 
\itemitem{\rm(i)} if $|H_\Omega|\le|E|<|\Omega|$, then a perimeter
minimizer
$E$ is the
union of all balls in $\overline \Omega$ of curvature equal to the curvature of
$\partial E\cap\Omega$, 
\itemitem{\rm(ii)} if $|B_\Omega|<|E|<|H_\Omega|$, then $E$ is the
union of two largest balls in $\Omega$, 
\itemitem{\rm(iii)} if $0<|E|\le |B_\Omega| $,
then $E$ is a ball.} 
\medskip
{\bf Proof.} 
Smooth $\partial \Omega$ as before but without requiring
the great circle condition. The same
regularity properties hold as before for  perimeter minimizers $E_n$
in the smoothed domains $\Omega_n$. Note that there is no singular
set since $n=2$. Also $\partial E_n\cap \Omega$ consists of circular
arcs. Thus if $x\in \Omega$  is a limit point of points
$x_n\in\partial E_n$ then it is easy to see geometrically that the
curvatures of the circular arcs in  $\partial E_n\cap \Omega$ must be
uniformly bounded in $n$. Regularity and curvature results for  the
limiting perimeter minimizer follows as before.

We claim that any perimeter minimizer $E$ must be convex. First note
that  $E$ cannot have an infinite number of components since otherwise
$\partial E$ would contain a limit point of points in the boundaries
of distinct components of $E$ which would violate the regularity of
$\partial E$. In addition each component must be simply connected
since otherwise one could add a bounded component of the complement of
$E$ to $E$ which would reduce the perimeter of $E$ and increase its
measure. Thus a scaling argument as in the  proof of Theorem
\thmref{econvex} would violate the fact that $E$ is a perimeter
minimizer. 

Also each component must be convex. To see this note that locally
$\partial E$ is a graph of a $C^{1,1}$ function $f$. Thus $f'$ is
Lipschitz continuous, monotone increasing (if axes are chosen
properly) on $f^{-1}(\partial E\cap\partial\Omega)$, and monotone
increasing on each component of $f^{-1}(\partial E\cap \Omega)$ from
which the claim easily follows.

Finally given two components  considering the two unique lines
which are support lines for both components one sees that one of the
components can be translated with out leaving $\Omega$ until it first
touches another component. This translation does not change the
measure of the overall set and does not increase perimeter so a new
perimeter minimizer is created. Due to the regularity of $\partial E$
the point of contact lies in $\Omega$. However this contradicts the
fact that the boundary of a perimeter minimizer must be a circular arc
locally in $\Omega$. Consequently there must be only one component
which we have already shown to be convex so $E$ is convex as claimed.

To establish the uniqueness and nestedness properties it is
sufficient to characterize perimeter minimizers. In fact we claim that
if $E$  is a perimeter minimizer with $|H_\Omega|\le|E|<|\Omega|$
then it is the union of all balls in
$\overline\Omega$ of curvature given by the curvature of $\partial
E\cap\Omega$. We prove the claim in two parts. We first establish that
if a point $x$ lies in $E$ then $x$ lies in a ball
contained in $\overline\Omega$ whose boundary has
the same curvature as $\partial E\cap\Omega$ . We finish by proving
that if $|H_\Omega|\le|E|<|\Omega|$ then $E$ contains all balls  with
the same curvature as $\partial E\cap\Omega$.

Assume $x\in E$ and let $d=\hbox{dist}(x,\,\partial E),\, r={1\over
k}$  where $k$ is the curvature of  $\partial E\cap\Omega$. If $d\ge
r$ then $x$ is clearly in a ball of radius $r$  contained in $\overline\Omega$
as claimed. If $d<r$ then choose a point $y\in\partial E$ closest to
$x$. Choose axes so that $y$ is the origin, $x$ lies on the positive
horizontal axis, and the vertical axis is tangent to $\partial
E$ at $y$. Let $C$ be the upper half of the circle of radius $r$
containing $y$  with center on the positive horizontal axis. Let
$(0,\,a)$ be the largest subinterval of $(0,\,2r)$
 over which the part of $\partial E$
lying above the horizontal axis is a graph. Let $f:(0,a)\rightarrow
\R$ be the function  having such a graph. Let
$g:[0,\,2r]:\rightarrow\R$ be the  function with $C$ as its graph.
Integrating the divergence form for curvature over $(\varepsilon,\,t)$
for $t<a,\,\varepsilon>0$ then letting $\varepsilon\rightarrow 0$
one obtains
$$
1-J(f'(t))=-\int^t_0 [J(f'(s))]'\,ds
\le \int^t_0 k\,ds
= -\int^t_0 [J(g'(s))]'\,ds
=1-J(g'(t))
$$
where $J(x)=x/(\sqrt{1+x^2})$ since $f'(0)=g'(0)=\infty$ (recall $\partial E$
is $C^{1,1}$).  Since $J(x)$ is monotone increasing this implies that
$g'(t)\le f'(t)$ on $(0,\,a)$. However $0=g(0)\le \lim_{s\to 0}f(s)$ so
$g(t)\le f(t)$ on $(0,\,a)$. From the estimate on $f'$ and the convexity of
$E$ we see that $a=2r$.  A similar argument shows that the part of $\partial E$
lying below the horizontal axis in fact lies below the other half of the 
circle of radius $r$ mentioned above. Thus from the convexity of $E$
we see that this circle lies in $E$ as claimed. Consequently
$E$ lies in the union of all balls of radius $r$ which lie in $\overline\Omega$.

To prove our second claim let $B$ be a ball of radius $r$ contained in
$E$ (such a ball exists by the above argument). 
Let $D$ be any other ball of radius $r$
contained in $\Omega$ and let $H$ be the convex hull of $B,\, D$.
Assume that $D$ is not a subset of $E$ so there exists $x\in
D\setminus E$. Thus $\partial E$ separates $x$ from $B$. However
$H^\circ$, the interior of $H$, lies in $\Omega$ so $\partial E\cap
H^\circ$ is locally a circular arc of radius $r$. The only way a
circular arc of radius $r$ can separate $x$ from $B$ is if it is a
half circle $C$ tangent at its end points to the line segments in
$\partial H$. In such a case $\partial E$ must contain $C$ and the
(possibly empty) line segments in $\partial H$ with endpoints in $C$
and  $\partial B$. Since $x\notin E$ one can translate $E$ towards $x$
 while remaining in $\Omega$ due to the geometrical relationship 
 between $E$ and $H$. The translated set is thus still a perimeter
 minimizer with end opposite $D$ lying in $\Omega$. Thus the end
 opposite $D$ is a circular arc and $E$ is the convex hull of two
 (possibly identical) balls of radius $r$. 

If $|E|\le |B_\Omega|$ then clearly $E$ is a ball. If
$|B_\Omega|<|E|\le|H_\Omega|$ then $\Omega$ satisfies a great circle
condition since the line segments in $\partial H$ must lie in $\partial\Omega$.
Thus we can use the characterization of
$E$ in Theorem \thmref{notstrict} as the convex hull of two largest balls in 
$\overline\Omega$. If $|H_\Omega|<|E|$ then $E$ cannot be a ball or a hull of
two balls in $\Omega$ as concluded in the last paragraph. Consequently
the assumption that $D$ was not a subset of $E$ is false and $D\subset
E$. Since $D$ was an arbitrary ball of radius $r$ we see that the
union of all such balls lies in $E$. Combining this with our earlier
conclusion we see that $E$ is in fact equal to the union of all such
balls.

The disjointness property for boundaries of perimeter minimizers
follows from nestedness of minimizers and the fact that  the curvature
of the boundary of a perimeter minimizer in $\overline\Omega$ strictly
increases as a function of the measure $v$ of the minimizer
if $|H_\Omega|\le v$, a fact which
follows directly from the characterization of minimizers. If  $E,\, F$
are minimizers with $E\subset F$, and $\partial E\cap\partial F
\cap\Omega$ is not empty then geometrically the curvature of $\partial
F\cap\Omega$ can not be larger than the curvature of $\partial
E\cap\Omega$. However this contradicts the monotonicity of curvature
as a function of measure mentioned above. 

\bigskip

\bigskip

\sectionnumber=4
\theoremnumber=0
\equationnumber=0
{\bf 4. Eqimeasurable Convex Rearrangement}
\bigskip
Various standard symmetrizations have the useful property of
rearranging  functions in an equimeasurable fashion while reducing
various norms such as $\|u\|_{L^p(\Omega)}+\|\nabla u\|_{L^p(\Omega)}$
and $\|u\|_{BV(\Omega)}$  (see \eqnref{np.1}).  However they alter
$\Omega$, the  domain of definition of $u$, unless  $\Omega$ has
appropriate symmetries. This is unfortunate from the point of view of
studying minimizers to  certain variational problems. Using results of
Section 2 we introduce an equimeasurable rearrangement which
preserves convex domains, reduces $\|u\|_{BV(\Omega)}$, and creates
level sets which are boundaries of convex sets, when $u\in BV(\R^n)$
with $u\ge0$ and $u=0$ in $\R\setminus\Omega$. Results of [LS]
imply
that such a rearrangement cannot exist for the norm
$\|u\|_{L^p(\Omega)}+\|\nabla u\|_{L^p(\Omega)}$, $p>1$. Any
equimeasurable rearrangement clearly fixes the first term in the BV
norm \eqnref{np.1}. From the co-area formula we will see that a
rearrangement which minimizes the perimeter of sets $\{u>t\}$  will
minimize the BV norm over an appropriate class of equimeasurable
functions . 

In minimizing functionals such as 
$$
 \|u\|_{BV(\Omega)} +\int_\Omega F(u) 
 + \int_0^{|\Omega|} G(u^*,{u^*}')\eqnlbl{norm2}
$$
over appropriate function classes , where $u^*$ is the decreasing
rearrangement of $u$, $u^*(v)=\sup\{t\,:\,|\{u>t\}|\ge v\}$, it is sometimes
straight forward  to derive regularity estimates for $u^*$. Assuming
continuity of $u^*$ the results of
Theorem \thmref{rearrangement} imply continuity for 
minimizers of \eqnref{norm2} in $\Omega\setminus H_\Omega$,
using the continuity and uniqueness properties of $\tilde u$.
Of course to apply Theorem \thmref{rearrangement} it is necessary
that $u=0$ on $\Omega$ is a boundary condition for the variational
problem and that one can establish $u\ge 0$ in $\Omega$ for minimizers
for instance by using a truncation argument.
Behaviour in $H_\Omega$ is also highly constrained by the characterization
of level sets up to translation. It is fairly straight forward
but more delicate to prove partial regularity results for $\nabla u$
if $\Omega\subset\R^2$ by analyzing interactions between boundaries of 
perimeter minimizers and $\partial \Omega$. However in higher dimensions
this is a difficult open problem.

Assume that $\Omega$ is a bounded convex set in $\R^n$. In addition
assume that $n=2$, or  $\Omega$
satisfies a great circle condition. Thus from Section 2 we have a
family of convex nested perimeter minimizers $E(v)$ defined as
follows. If $B_\Omega$ is a largest ball in $\Omega$ and $H_\Omega$
is the union of all such balls then if $0<v\le|B_\Omega|$ let $E(v)$
be a ball  of measure $v$ centered symmetrically in $H_\Omega$, if
$|B_\Omega|<v\le|H_\Omega|$ (in which case $\Omega$ satisfies a great
circle condition) then let $E(v)$ be  the convex hull of two largest
balls symmetrically centered in $H_\Omega$ and of measure $v$, finally
if $|H_\Omega|<v<|\Omega|$ let $E(v)$ be the unique perimeter
minimizer of measure $v$ shown to exist in Section 2.

Define
$$
BV^+_0(\Omega)=
\{u\in BV(\R^n)\,:\,u\ge0,\, u=0 \hbox{ in } \R^n\setminus\Omega\}
$$
and define the convex rearrangement of a function $u\in BV^+_0(\Omega)$ by
$$
\tilde u(x) = \inf\{s\ge 0\,:\, x\notin E(|\{u>s\}|)\}.
$$
\bigskip
{\bf\thmlbl{rearrangement} Theorem.} {\it 
If $\Omega$ is as above and $u\in BV^+_0(\Omega)$
then $\tilde u$ is upper semicontinuous in $\R^n$, continuous in
$\Omega$ if $u^*$ is continuous 
(equivalently $|\{u>t\}|$ is strictly increasing),
$\tilde u\in BV^+_0(\Omega)$,
$$
|\{\tilde u>t\}|=|\{u>t\}|
$$
for all $t$, and
$$
\norm{\tilde u}_{BV(\R^n)}\le \norm{u}_{BV(\R^n)}.\eqnlbl{bvinequality}
$$
If there is equality in \eqnref{bvinequality} then $\tilde u= u$ 
in $\R^n\setminus H_\Omega$ in the BV sense, and in $H_\Omega$
the  sets $\partial^*\{\tilde u>t\}$ are translations of
$\partial\{u>t\}$.}
\smallskip
{\bf Remark:} From the remark after Theorem \thmref{enested}
one sees that it is possible to create a rearrangement even if
the convexity assumption is relaxed. However it is unclear that
one can in this context establish qualitative information 
analogous to convexity of $\{\tilde u>t\}$.

\medskip

{\bf Proof.} Semicontinuity and continuity results are clear from the
definition of $\tilde u$ and the disjointness results on boundaries
of perimeter minimizers in $\Omega$.
It is also clear that $\tilde u\in BV^+_0(\Omega)$.
Due to the convexity and  nestedness (which is strict in $\Omega$)
of the sets $E(v)$ we see that

$$
E^\circ(|\{u>t\}|)\subset \{\tilde u>t\}\subset E(|\{u>t\}|)
$$
thus 
$$
|\{u>t\}|=|E(|\{u>t\}|)|= |\{\tilde u>t\}|.
$$
and 
$$
P(\{\tilde u>t\})=P(E(|\{u>t\}|))\le P(\{u>t\}).
$$
The result on BV norms then follows from the co-area formula.

If one has equality in the BV norm expression then from
the co-area formula and the minimization property
of the sets $E(v)$ it is clear that 
$P(\{\tilde u>t\})= P(\{u>t\})$, and consequently $\{\tilde u>t\}$ 
is a perimeter minimizer for almost all $t$.
Let $t_0=\sup\{t\,:\, |\{\tilde u>t\}|\ge |H_\Omega|\}$ so 
applying the uniqueness result for perimeter minimizers
we see that $\{\tilde u>t\}$ and $\{u>t\}$ have the same
measure theoretic closure for almost every $t,\, 0\le t <t_0$. For $t\ge t_0$
we have $|\{\tilde u>t\}|<|H_\Omega|$ so
this is true only up to translation within $H_\Omega$
in which case $\partial^*\{\tilde u>t\}$ is a translation of $\partial\{u>t\}$
(recall $\{u>t\}$ is convex) as claimed. 
This is easily justified for all $t,\,t\ge t_0$
by a limit argument.

Returning to the case $0\le t<t_0$ let 
$E$ be an arbitrary measurable subset of $\Omega$
and $d\mu=\upchi_E\, dx$ where $dx$ represents Lebesgue measure. 
From Fubini's theorem
we see that
$$
\int_0^{t_0}\mu(\{u>t\})\, dt=\int\int_0^{t_0}\upchi_{\{u>t\}}\,dt\,d\mu
=\int_E \min(u,\,t_0).
$$
Using the fact that  $\{\tilde u>t\}$ and $\{u>t\}$ have the same
measure theoretic closure for almost every $t,\, 0\le t <t_0$
we conclude that $\min(u,\,t_0)=\min(\tilde u,\,t_0)$ almost
everywhere. Recalling that $\{u>t\}$ and the set theoretic closure of $\{u>t\}$ 
are subsets of $H_\Omega$ for $t>t_0$ it 
is clear that $\tilde u=u$ almost everywhere
in $\R\setminus H_\Omega$.\qed
\vfill\eject
\centerline{\bf References}
\bigskip
\parindent50truept
\flushitem{[BK]} Brezis, H., and D. Kinderlehrer, {\it The
smoothness of solutions to nonlinear variational inequalities}, 
Ind. Univ. Math. J., 23(1974), 831 -- 844.
\bigskip
\flushitem{[Fe]} Federer. H., {\it Geometric measure theory},
   Springer Verlag, New York, Heidelberg, 1969.
\bigskip
\flushitem{[GMT1]} Gonzalez, E., U. Massari and I. Tamanini, {\it
Minimal boundaries enclosing a given volume}, Manuscripta Math.,
34(1981), 381 -- 395.
\bigskip
\flushitem{[GMT2]} Gonzalez, E., U. Massari and I. Tamanini, 
{\it On the regularity of sets minimizing perimeter with a volume
constraint}, Ind. Univ. Math. J., 32(1983), 25 -- 37.
\bigskip
\flushitem{[Gr]} Gr\"uter, M., {\it Boundary regularity for
solutions of a partitioning problem}, Arch. Rat. Mech. Anal., 
97(1987), 261-270.
\bigskip
\flushitem{[Gi]} Giusti, E., {\it Minimal surfaces and functions of
bounded variation}, Birkh\"auser, 1985.
\bigskip
\flushitem{[LS]} Laurence, P., E.W. Stredulinsky, {\it On Quasiconvex
Equimeasurable Rearrangement, a Counterexample and an Example},
J. fur Reine Angew. Math., 447 (1994), 63-81.
\bigskip
\flushitem{[MM]} Massari, U. and M. Miranda, {\it Minimal surfaces of
codimension one}, Mathematics Studies, North Holland, 91 (1984).
\bigskip
\flushitem{[S]} Simon, L., {\it Lectures on geometric measure
theory}, Proc. Centre Math. Analysis, ANU, 3 (1983).
\bigskip
\flushitem{[T]} Tamanini, I., {\it Boundaries of Caccioppoli sets
with H\"older-continuous normal vector}, J. fur Reine Angew. Math.,
334 (1982), 27-39.
\bigskip
\flushitem{[Z]} Ziemer, W.P., {\it Weakly differentiable functions,}
Springer-Verlag, GTM Series, 120 1989.

\bye